\documentclass{iopart}
\usepackage{iopams}

\usepackage{geometry}
\usepackage{epstopdf}
\usepackage{latexsym}
\usepackage[pdftex]{graphicx}
\usepackage{epstopdf}
\eqnobysec

\begin{document}

\title[]{Bifurcation from a normally degenerate manifold}

\author{D.R.J.~Chillingworth}
 
\address{School of Mathematics, University of Southampton,
       Southampton SO17 1BJ
        }
 
\author{L.Sbano}

\address{Mathematics Institute, University
 of Warwick, Coventry CV4 7AL}


\begin{abstract}
Local bifurcation theory typically deals with the response of a
degenerate but isolated equilibrium state or periodic orbit 
of a dynamical system to
perturbations controlled by one or more independent parameters, and
characteristically uses tools from singularity theory.  There are many
situations, however, in which the equilibrium state or periodic orbit
is not isolated but belongs to a manifold $S$ of such states, 
typically as a result of continuous symmetries in the problem.  
In this case the bifurcation analysis requires a
combination of local and global methods, and is most tractable in the
case of normal nondegeneracy, that is when the degeneracy is only
along $S$ itself and the system is nondegenerate in directions normal to $S$.  

In this paper we consider the consequences
of relaxing normal nondegeneracy, which can generically occur within
$1$-parameter families of such systems.  We pay particular attention to
the simplest but important case where $\dim S=1$ and where the normal
degeneracy occurs with corank $1$. Our main focus is on uniform
degeneracy along $S$, although we also consider aspects of the
branching structure for solutions when the degeneracy varies at
different places on $S$. The tools are those of singularity theory
adapted to global topology of $S$, which allow us to explain the
bifurcation geometry in natural way. In particular, we extend and
give a clear geometric setting for earlier analytical results of 
Hale and Taboas.
\end{abstract}
\maketitle


\def\rem{\noindent{\sl Remark.\quad}}
\def\rems{\noindent{\sl Remarks.\quad}}

\def\qed{\hfill$\Box$}

\def\ssk{\smallskip}
\def\msk{\medskip}
\def\bsk{\bigskip}
\def\noi{\noindent}
\def\br{{\mathbb{R}}}
\def\bc{{\mathbb{C}}}
\def\bn{{\mathbb{N}}}
\def\bz{{\mathbb{Z}}}
\def\bq{{\mathbb{Q}}}
\def\bs{{\mathbb{S}}}
\def\del{\delta}
\def\eps{\varepsilon}
\def\x{\xeta}
\def\om{\omega}
\def\ga{\gamma}
\def\sig{\sigma}
\def\pp{\partial}
\def\codim{\mbox{codim}}
\def\A{{\mathcal A}} 
\def\H{{\mathcal H}} 
\def\K{{\mathcal K}} 
\def\M{{\mathcal M}} 
\def\V{{\mathcal V}} 
\def\J{{\mathcal J}}
\def\F{{\mathcal F}}  
\def\R{{\mathcal R}}  
\def\tilF{\widetilde{{\mathcal F}}}  
\def\ba{\bar a}
\def\bb{\bar b}
\def\noi{\noindent}

\def\ots{0\times S}

\def\rmbb{{\textrm(}}
\def\rmbe{{\textrm)}}

\def\spn{\mathop{\mathrm{span}}\nolimits}
\def\grad{\mathop{\mathrm{grad}}\nolimits}
\def\rank{\mathop{\mathrm{rank}}\nolimits}
\def\dbd#1#2{\frac{\partial{#1}}{\partial{#2}}}

\def\proof{\emph{Proof.\enspace}}
\def\endproof{\hfill\qed}

\def\qed{\hfill$\Box$}

\def\nhd{neighbourhood\ }
\newcommand\til{\tilde}

\def\ker{\mathop{\mathrm{ker}}\nolimits}
\def\Lin{\mathop{\mathrm{Lin}}\nolimits}
  
\newtheorem{lem}{Lemma}[section]
\newtheorem{prop}[lem]{Proposition}
\newtheorem{cor}[lem]{Corollary}
\newtheorem{theo}[lem]{Theorem}
\newtheorem{conj}[lem]{Conjecture}
\newtheorem{defn}[lem]{Definition}
\newtheorem{obs}[lem]{Observation}

\section{The general problem}\label{sect:gen}
\subsection{Introduction} \label{ss:intro}
Let $F_0:M\to TM$ be a smooth\footnote{By smooth we mean $C^\infty\,$.  However,
it is clear that many of our results hold with only a few orders of
differentiability.} vector field on a smooth (possibly
infinite-dimensional) manifold $M$, with
the property that the flow generated by $F_0$ has a smooth 
$d$-dimensional manifold $S$ of equilibria with finite
$d>0$. In other words, $F_0$ vanishes everywhere on $S$.  
We now perturb $F_0$ to 
\begin{equation} 
\label{e:orig}
F_\eps:=F_0 + \eps \F + O(|\eps|^2)
\end{equation}
where $\F$ is a smooth vector field on a \nhd of $S$ in $M$ and where
$\eps$ is a real parameter, and we
ask what can be said about zeros of $F_\eps$ close to $S$ in $M$ for
$\eps$ sufficiently small but nonzero.  In particular, we ask which
points $x_0\in S$ are such that a branch of equilibria emanates from
them as $\eps$ varies away from zero.
\ssk

More generally we may consider a multiparameter deformation $F_\eps$
where $\eps\in\br^q$, as in \cite{HA2}\cite{V1} or \cite{C2,C4} for example. 
In this case the term $\eps \F$ in (\ref{e:orig}) is interpreted as
$\sum_{i=1}^q\eps_i\F_i$ for a collection $\{\F_1,\ldots,\F_q\}$ of smooth
vector fields, and $O(|\eps|^2)$ means terms of degree 2 or higher in the
components of $\eps$. Here it is useful to think in terms of
equilibrium branches as $\eps$ moves away from the origin \emph{in a given
direction} in parameter space $\br^q$, and we discuss this more carefully
below.  Of course, it is not clear \emph{a priori} that branching
behaviour depends only on $\F$ and not on the $O(|\eps|^2)$ terms,
although it is a key aspect of our results
that often under suitable generic\footnote{In our context a generic
  condition or assumption is one which holds for an open dense set of
  the maps or objects under consideration.  Making this precise
  requires formal machinery of transversality in jet spaces and
  techicalities which we mostly prefer to avoid: in our cases the
  interpretation of genericity is straightforward.}      
 assumptions on $\F$ this is indeed the case.
\ssk

In many situations where bifurcation from a manifold 
arises in applications the manifold $S$ is
{\em normally non-degenerate} with respect to $F_0\,$, meaning that at each 
point $x\in S$ the tangent space $T_xS$ (which automatically lies in 
$\ker T_xF_0$) in fact coincides with $\ker T_xF_0$, 
so that $T_xF_0$ is as non-degenerate as
possible under the given geometric assumptions.  Bifurcation under
these conditions has been well studied:
see~\cite{Amb},\cite{C4},\cite{CC3},\cite{D1},\cite{HA2},\cite{V1}.
A standard  technique is to show by a global Lyapunov-Schmidt reduction along $S$
that for small $\eps$ there is a smooth manifold $S_\eps$ close to $S$ 
(viewed as the image of a small section of the normal bundle $NS$ of $S$ in
$M$) such that the `normal' component of $F_\eps$ vanishes on $S_\eps$.  
This reduces the problem to a study of the `tangential' component of $F_\eps$,
which corresponds to a tangent vector field $\til F_\eps$ on $S$
via the projection $NS\to S$.  
The sought-for zeros of $F_\eps$ close to $S$ then correspond to
the zeros of the vector field $\til F_\eps$ on $S$
itself.\footnote{For a more careful description see Section~\ref{sect:LS}.}

In the absence of further information the number of zeros of 
$\til F_\eps$ can often be given a lower bound by topological
means such as index theory: see the techniques of \cite{D1,D2} for
example.  In the variational case where $F_\eps$ 
is the gradient field of a real-valued function the same may be
assumed for $\til F_\eps$ (see Section \ref{ss:gradient}) and
the usual tools are Lyusternik-Schnirelmann category theory in general or
Morse theory if critical points are assumed to be nondegenerate: see
for example \cite{Amb},\cite{Chang},\cite{RE2},\cite{D2},\cite{CHH} and
\cite{WE},\cite{ACE,ACZ} for applications to Hamiltonian systems.
With further information on symmetries \cite{D3},\cite{MO},\cite{B2},\cite{V1} 
or on the actual perturbing field $\F$ 
more specific results can be obtained as in
\cite{C1,C2,C3,C4} and \cite{CC1,CC2,CC3}.
\msk

The purpose of this paper is to investigate the consequences of 
relaxing the normal non-degeneracy condition by supposing that $T_xS$ 
is a proper subspace of $\ker T_xF_0\,$, that is
\begin{equation} \label{e:ccorank}
\dim\ker T_xF_0 = k + d
\end{equation}
for all $x\in S$, where $k\ge1$ does not depend on $x\,$.  
We call this {\em normal degeneracy} with {\em constant corank}  $k$. 
Our most detailed results are for the simplest (but important) case
when $k=1$, but some of the methods and observations are applicable
also to the general case.
\msk

\subsubsection{The group orbit case.}
There is a natural situation in which manifolds of equilibria arise automatically, 
namely where the manifold $S$ is the orbit of the action of a Lie group 
$\Gamma$ under which $F_0$ is equivariant. 
In particular this occurs when $F_0$ is the 
gradient of a $\Gamma$-invariant function $f_0:M\to\br$.
In problems formulated on function spaces of periodic functions, 
there is a natural circular symmetry defined by the phase-shift, and
in mechanical systems natural symmetries arise from coordinate
invariance. Translational symmetry can also play a key role in
bifurcation theory of differential equations on the real line: 
see~\cite{Mag} for example. Problems that we consider, in which the
perturbed system $F_\eps$ does not share the symmetry ($\,\Gamma$-equivariance)  
of $F_0\,$, go under the general heading of {\em forced symmetry-breaking}.
 
For a typical group orbit we expect normal non-degeneracy,
but if $F_0$ is a member of a smooth $1$-parameter
family of vector fields $F^\gamma_0$ each of which is
$\Gamma$-equivariant then, since (for square matrices)
dropping rank by $1$ is a codimension-$1$ phenomenon, we would expect to find
isolated values of the parameter $\gamma\in\br$ for which the corresponding
equilibrium manifold $S^\gamma$ is normally
degenerate with constant corank $1\,$.  The constancy of corank 
in this case is a consequence of the equivariance.
\msk

\subsubsection{Non-constant corank.}
In contrast with the above, 
non-degeneracy may fail in that
$\ \dim\ker T_xF_0 > \dim S\ $ only at isolated points $x\in S\,$.
This does not arise in the group orbit case, since the group
action implies that the behaviour at all points of $S$ is `the same'.
Analysis of this type of degeneracy requires different techniques, which
we touch upon in Section~\ref{s:q=k=1}.  See also \cite{RC}.
\subsubsection{Sources of the problem.}   
The break-up of families of periodic orbits of Hamiltonian systems is a
classical problem, first treated in a modern spirit in the fundamental
papers of Weinstein~\cite{WE} and Moser~\cite{MO}.  
Abstract treatment of bifurcation from a manifold was pioneered by
Hale~\cite{HA2}, with applications to planar forced ordinary
differential equations, and extended by Hale and Taboas~\cite{HT2} to
the case of normal degeneracy.  Further important techniques and
results relating particularly to symmetries were set out by 
Dancer~\cite{D1,D2,D3} and Vanderbauwhede~\cite{V1}, again supposing
normal nondegeneracy.  Applications of bifurcation from a manifold
were exploited by Ambrosetti \emph{et al.}~\cite{ACE,ACZ} in the
context of Hamiltonian systems and also more generally~\cite{Amb}.  
Meanwhile, a general geometric framework for multiparameter
bifurcation from a normally nondegenerate manifold was
formulated by the first author~\cite{C1,C2,C4}, and customised tools
adapted to families of periodic orbits of ordinary differential
equations were systematically refined in a series of papers by
Chicone~\cite{CC1,CC2,CC3}. It is the work of Hale \emph{et al.} and
Chicone that is the main inspiration for the current paper.
%
\subsection{Outline of the approach}
An overview of our approach in this paper is as follows.
First we apply Lyapunov-Schmidt reduction to enable us to focus on the
core of the problem.  As a consequence of normal
degeneracy, the reduced problem still involves a $k$-dimensional
normal component of the vector field, as well as the automatic
$d$-dimensional tangential component.
We would like to `unfold' this degenerate $(k+d)$-dimensional vector field,
using $\eps$ as a multi-parameter, but we are unable to use
standard tools of unfolding (deformation)
theory directly as the degeneracy along
$S$ causes the local problem to have infinite codimension.  
Therefore we make appropriate generic assumptions about the behaviour
of the deformation term $\F$ along the manifold $S$.  With this,
we can in principle describe generic geometry of such deformations,
at least for small $k\,$.  In the fundamental
case $k=1$ we are able to give a rather complete picture.
\msk

In Section~\ref{sect:LS} we formulate the initial reduction process, and
in Section~\ref{s:versal} we show how to adapt the standard notion of versal 
deformation to our context, 
giving particular attention to the simplest nontrivial low-dimensional case.
In Section~\ref{s:d=k=1} we combine global geometry with local 
algebra to describe generic branching structure, and in
Section~\ref{s:q=k=1} we study local branching structure in closer detail.
Finally, in Section~\ref{s:apps} we discuss three areas of application
of the methods set out in this paper.
\section{Lyapunov-Schmidt reduction}
\label{sect:LS}
%
%
Let $NS$ be a normal bundle for $S$ in $M$, 
that is a vector sub-bundle of $T_SM:=TM|_S$ orthogonal to $TS$ (using 
a Riemannian structure $\nu$ on $M$). 
A \nhd of the zero section of $NS$ can be taken to parametrise 
a tubular \nhd of $S$ in $M\,$ via the exponential map associated to
$\nu$,  and so then we may without loss of generality
regard the whole bifurcation problem
as taking place on a \nhd $W$ of the zero section of the bundle $NS$.
\msk

Let $H$ be a connection on $NS$, that is a (smooth) sub-bundle of
$T\,NS$ orthogonal to (vectors tangent to) the fibres of $NS$: thus
\[
T_wNS=H_w\oplus T_wN_x=H_w\oplus N_x
\]
for all $w$ lying in the fibre $N_x$ of $NS$ over $x\in S$. This
enables us to define a linear bundle map $\pi:T\,NS\to NS$ as
projection onto the second factor in each tangent subspace $T_wNS$.
\msk

Given a vector field $F:W\subset NS\to TNS$ the composition
$\pi F:W\to NS$ satisfies $\pi F(W\cap N_x)\subset N_x$ for
each $x\in S$.  Thus $\pi F$ is a `vertical' vector field on $W$
obtained by projecting $F$ parallel to the `horizontal' subspaces
determined by $H$.  
If $w\in W\cap N_x$ we write $DF(w)\in Lin(N_x,N_x)$ to denote the
`vertical' derivative of $F$ (i.e. the derivative of
$F|{N_x}$) at $w$.  
\msk

Let $F_{\eps}:M\to TM\,,\,\eps\in\br^q$ be a vector field
as in Section~\ref{sect:gen}, with $F_0(x)=0$  for all $x\in S$.
As above, we replace $M$ by $W\subset NS$.
Since $T_xS\subset \ker T_xF_0$ for
all $x\in S$, the constant corank hypothesis~(\ref{e:ccorank}) implies
that at the origin $0_x$ of $N_x$ the kernel $K_x:=\ker DF_0(0_x)$
is a $k$-dimensional linear subspace of $N_x$ for every $x\in S\,$.
We now make the following assumptions:
\bsk

\begin{itemize}
\item[(A1)] The subspaces $\{K_x\}$ together 
define a smooth $k$-dimensional vector sub-bundle $K$ of $NS$.
\item[(A2)]  There is a smooth codimension-$k$ vector sub-bundle
$L=\{L_x\}$ of $NS$ complementary to $K$ in $NS\,$; 
thus $N_x=K_x\oplus L_x$ for all $x\in S$. 
\item[(A3)]  There is a smooth vector bundle decomposition
$NS=P\oplus R$ such that $R_x$ is the range of $DF_0(0_x)$
and $\dim P_x=k$ for all $x\in S$.
\end{itemize}
\bsk

\rem
In finite dimensions the smoothness of $K$ holds automatically, given the
constant dimension $k$, since the kernel of a matrix is the
orthogonal complement of the range of its transpose: thus if $DF_0(0_x)$
varies smoothly with $x$ then a basis for $\br^n$ orthogonal to the
rows of $DF_0(0_x)$ can also be chosen to vary smoothly.  The construction
of $L$ and $P$ is straightforward.

Smoothness of $K$ also holds in many naturally-arising 
infinite-dimensional cases where $DF_0(0_x)$ is a differential operator,
so that its kernel consists of the solutions $\xi(t)$ to a
$k$-dimensional system of linear differential equations; 
construction of $L$ and $P$ is also often easy in such cases.

If $DF_0(0_x)$ is self-adjoint (as for a gradient vector field, for
example) then we may take $L_x=R_x$ and $P_x=K_x$.
\bsk

\begin{figure}[htbp] 
   \centering
   \includegraphics[scale=0.5]{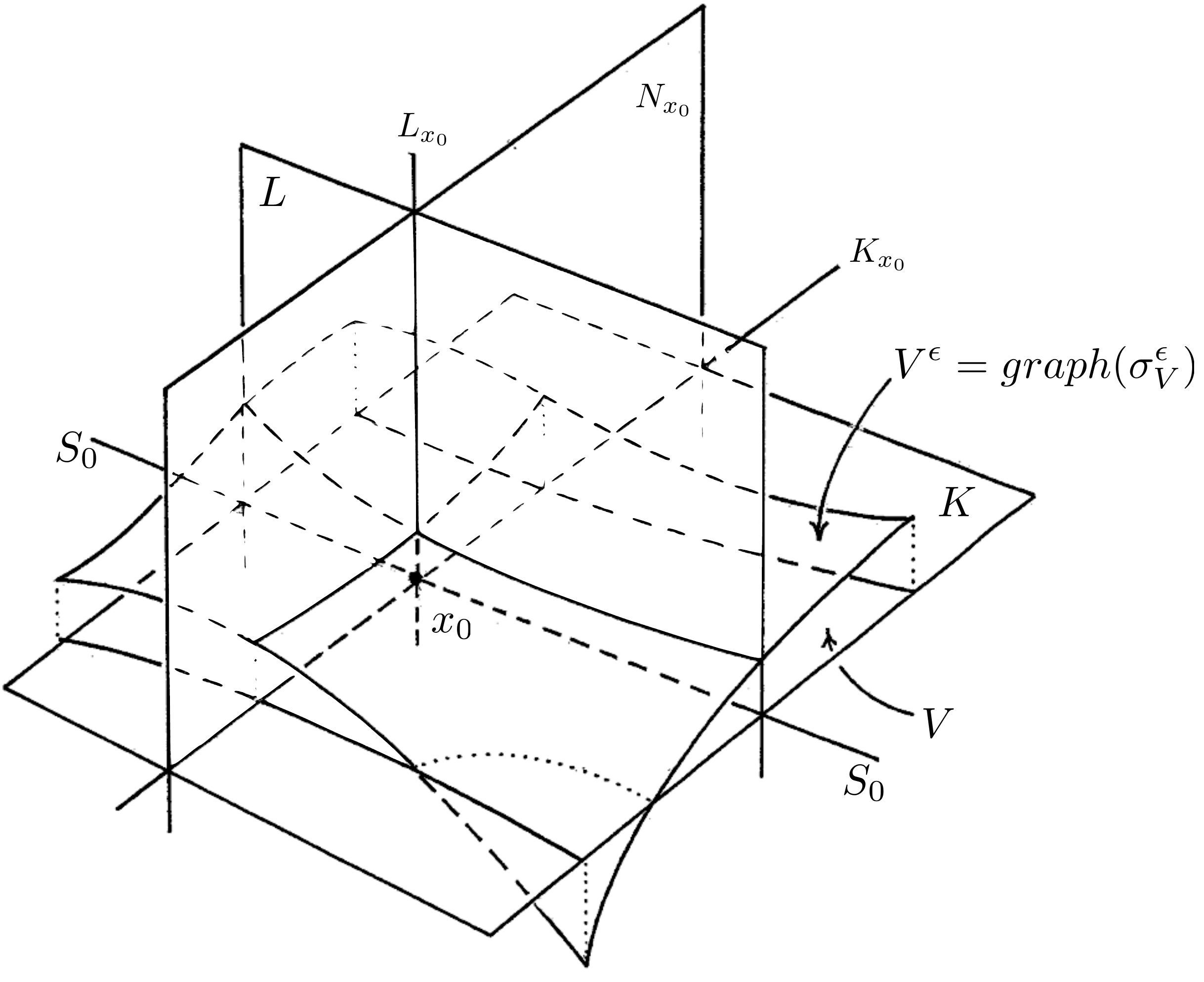} 
  \caption{The geometry of Lyapunov-Schmidt reduction at $x_0\in
 S_0\,$. Here $N$ is the normal bundle to $S=S_0$ in $M$ and 
 $K_{x_0}=\ker DF(x_0)$.}
   \label{FigureAA}
\end{figure}

Let $S_0$ denote the (image of the) zero section of the bundle $K$. 
Of course $S_0=S$, but the $S_0$ notation will 
focus attention on $S$ as a submanifold of $K$.
%
%
%
Applying Lyapunov-Schmidt reduction 
globally along a \nhd of $S_0$ in $K$
we obtain the following result.  
\begin{prop}    \label{p:lsr}
There exists a \nhd\ $U$ of $S_0$ in $K$ with the property that
for all sufficiently small $\eps$ there is a unique smooth bundle map
$\sigma^\eps:U\to L$ such that 
\begin{enumerate}
\item $\sigma^0|S_0=0$ and
\item the $R$-component of $F_\eps$ vanishes precisely at the points lying on
the graph $U^\eps\subset NS$ of $\sigma^\eps$.  
\end{enumerate}
\end{prop}
\proof
Let $r$ denote the bundle projection $NS=P\oplus R\to R$, 
and consider the composition $r\circ \pi F_\eps:W\to R$.  
When $\eps=0$ the vertical derivative of this smooth map at $0_x\in N_x$
(identified with $x\in S_0$) is by definition of $R$
a linear surjection with finite-dimensional
kernel $K_x$.  Hence by the Implicit Function Theorem (IFT),
given $x_0\in S_0$ there is a neighbourhood $V$ of $x_0$ in $K$ 
and $\eps_0>0$ such that for all $|\eps|<\eps_0$ and $x\in V\cap S_0$
the solution set to $r\circ \pi F_\eps=0$ 
in the fibre $N_x$ is the intersection with $N_x$ of the
graph $V^\eps$ of a unique smooth bundle map $\sigma_V^\eps:V\to L\,$.  See 
Figure \ref{FigureAA}. The construction of $\sigma^\eps$
on a global \nhd $U$ of $S_0$ from local sections $\sigma_V^\eps$
then follows by uniqueness of $\sigma_V^\eps$ on $V$ and the
local compactness of $S$. 
\endproof
%
\medskip

Since the $(d+k)$-dimensional manifold $U^\eps$ is precisely where
the vertical component of the vector field $F_\eps$ has zero $R$-component, 
it follows that the zeros of $F_\eps$ are
the zeros of the vector field on $U^\eps$ obtained by projecting
$F_\eps$ to $TS\oplus P$. 
(Note that this projected vector field is not
necessarily tangent to $U^\eps$.)
The problem therefore distils  to that of solving the {\em reduced problem}
\begin{equation}  \label{e:red}
\til F_\eps(u):=({\bf 1}-r\circ\pi)F_\eps(u,\sigma^\eps(u))=0
                              \in TS\oplus P\,, \qquad u\in U\subset K\,
\end{equation}
where ${\bf 1}$ denotes the identity map on $TM=TNS$.  Here if
$w=(u,\sigma^\eps(u))\in N_x$ we identify $H_w$ with $T_xS$ through
the natural bundle projection $TNS\to TS$, so that ${\bf 1}-r\circ\pi$
is projection on $TS\oplus P$.
\msk

\rem   If $F$ is of class $C^p$ then $\til F$ is also $C^p$ for
$1\le p\le\infty,\omega$.
\msk

\subsection{The gradient case}     \label{ss:gradient}
In the important special case when $F_\eps$ is the gradient $\nabla f_\eps$ 
of a real-valued function $f_\eps$ (with respect to the Riemannian
structure $\nu$) then the zeros of 
$\til F_\eps$ correspond to the zeros of a gradient field on $U$.  
To see this, take $R=L$ and $P=K\,$ (see earlier remark)
and observe that if $u\in U$ and 
$\til f_\eps(u):=f_\eps(u,\sigma^\eps(u))$ then
\begin{equation}
\begin{array}{ll}
\nabla\til f_\eps(u)&=({\bf 1}-r\circ\pi)\nabla
                                            f_\eps(u,\sigma^\eps(u))\\
                    &=\til F_\eps(u),
\end{array}
\end{equation}
there being no contribution to $\nabla\til f_\eps(u)$
from the second argument of $f_\eps$ because $r\circ\pi \nabla f_\eps$ 
vanishes identically on $U^\eps$ by construction.
\subsection{Interpretation of the reduced problem} \label{ss:interpret}
The reduced vector field $\til F_\eps:U\subset K\to TS\oplus P$
satisfies
\begin{equation} 
     \til F_0(0_x)=0 \quad\mbox{ for all } x\in S \label{e:tilf}
\end{equation}
and
\begin{equation}
   D\til F_0(0_x)=0 \quad\mbox{  for all } x\in S, \label{e:dtilf}
\end{equation}
since $F_0=0$ and $\sigma^0=0$ on $S_0$, and also
$({\bf 1}-r\circ\pi)DF_0(0_x)=0$ by definition of $r$.
We can therefore regard $\til F_\eps$ as arising from a family of smooth singular
map germs
\[
\til F_0|K_x:K_x,0_x \to T_xS\oplus P_x
\]
parametrised (globally) by $x\in S$ which is then perturbed (locally) by
$\eps\in\br^q$ for $|\eps|$ small.
After shrinking the original \nhd $U$ if necessary, we may write
\begin{equation}  \label{e:calF}
\til F_\eps(u)=\til F_0(u)+ \eps\tilF(\eps,u) 
\end{equation}  
with $u\in U\subset K$ and $\eps\in\br^q$ by analogy
with the original problem in Section~\ref{sect:gen}. 
Note that $\til F_\eps(0_x)$ is not necessarily zero when $\eps\ne0$.
\subsection{Uniform degeneracy} \label{ss:uniform}
It may happen that the map germ $\til F_0|{K_x}$ is in some
precise sense independent of $x\in S$: we shall make this more explicit
in Section~\ref{s:versal} below. 
In this case we say the normal degeneracy is \emph{uniform}.  
For example, in the group orbit case with trivial isotropy (no
nontrivial element of $\Gamma$ keeps $x\in S$ fixed)
the action of $\Gamma$ provides a natural coordinate system for the whole of
$NS$ as a product bundle over $S$.  In these coordinates the decompositions
$N_x=K_x\oplus L_x=P_x\oplus R_x$ are independent of $x\in S$
and uniformity in a natural sense follows automatically. 
\medskip

Machinery to study perturbations of functions and maps
is provided by singularity theory.  In the next section we introduce 
the key ideas that will be required for our purposes.

\section{Versal deformations}  \label{s:versal}

For ease of exposition we make a simplifying assumption, namely
that all of the bundles over $S$ that are involved in the construction
of the reduced vector field $\til F_\eps$ (including the tangent
bundle $TS$ are trivial bundles.  This
is automatically the case if $S$ is contractible (a copy of $\br$ for
example) and is usually the case in applications when $S$ is a circle.
Even when the assumption fails, we can apply the methods below to
contractible open subsets of $S$ and piece together the results
subsequently.  We remark later on how the methods would need to be
adapted to apply intrinsically to nontrivial bundles.
\bigskip

Let us therefore take $K=S\times\br^k$ and 
$TS\oplus P=S\times\br^d\times\br^k$
so that as in Section~\ref{ss:interpret} we regard $\til F_\eps$ as a
family of smooth map germs
\[
\til F_\eps(x,{\cdot}):\br^k,0 \to \br^{d+k}
\]
parametrised by $x\in S$ and (small) $\eps\in\br^q$.
Since we are here interested in the \emph{zero locus}
\[
Z_{\til F}=\{(\eps,x,y)\in \br^q \times 
                 S\times \br^k:\til F_\eps(x,y)=0\in\br^{d+k}\}
\]
for small $|y|$ rather than other geometric features of 
$\til F_\eps$ as a map, we are free to
apply coordinate changes in a \nhd of $0\times S_0=0\times S\times 0$
that respect the geometry of the zero locus and the role of $\eps$ as
a parameter.
A natural class of coordinate transformations to use is the
class of \emph{contact equivalences}\footnote{Not to be confused with 
contact transformations from classical mechanics.} or 
$\K$-\emph{equivalences}.
Singularity theory provides powerful tools for handling classification
and geometry of map germs up to $\K$-equivalence, and the essence of
this paper is the adaptation of those techniques to our class of
bifurcation problem. It is in order to use the tools of singularity
theory that we assume $C^\infty$ smoothness of all data in the problem.
Background material for the singularity theory in this section
can be found in a range of texts such as~\cite{AGLV,AVG,GI} and from 
a slightly different viewpoint in~\cite{GSS}.\\
\subsection{Basic versal deformation}
We here introduce the main elements of the theory of $\K$-versal 
deformations of map germs that will be used in our analysis. 
To recall the concepts we shall use generic symbols $x$ and $F$ that
are not the same as $x\in S$ and the map $F_\eps$ of our main problem.
\begin{defn}
Two smooth map germs $f,g:\br^n,0\to\br^p,0$ are $\K$-{\em equivalent} 
if there is a diffeomorphism germ $\phi$ at the origin in $\br^n$ and an
invertible $p\times p$ matrix $A$ depending smoothly on $x\in\br^n$ 
such that\footnote{The symbol $\equiv$ indicates the two map germs
are identical: it is not an equation for $x$.}
\begin{equation}  \label{e:equiv}
g(x) \equiv  A(x) f(\phi(x)).
\end{equation}
\end{defn}
Clearly the zero set for $f$ is the image under $\phi$ of the zero set
for $g\,$.  Since $A$ is invertible any degeneracy of $f$ at the
origin is \lq the same' as that of $g\,$.
\begin{defn} \label{d:deform}
An $r$-parameter \emph{deformation} of a map germ $f:\br^n,0\to\br^p,0$ is
map germ $F:\br^r\times\br^n,(0,0)\to\br^p,0$  such that $F(0,x)\equiv f(x)$.
\end{defn}
To keep notation simple we shall write $F$ as a map, the fact that it
represents a map germ at $(0,0)$ being understood.
\msk

Now for the key definition through which we classify types of
bifurcation behaviour.  
\begin{defn}
Two deformations $F_1,F_2:\br^r\times\br^n\to\br^p$ of $f$
are $\K$-{\em equivalent} if there is a diffeomorphism 
germ $\Phi$ at the origin in $\br^r\times\br^n$ of the form
\begin{equation} \label{e:diffeo}
\Phi(a,x)=(\psi(a),\phi(a,x))
\end{equation}
and an invertible $\,p\times p$ matrix $A$ varying smoothly with respect to
$(a,x)\in \br^r\times \br^n$ such that
\begin{equation}  \label{e:Kequiv}
F_2(a,x)\equiv A(a,x)\,F_1\circ\Phi(a,x).
\end{equation}
\end{defn}
Note that $F_1,F_2$ are $\K$-equivalent as maps, but the coordinate change
$\Phi$ distinguishes between the roles of the variables $x$ and
parameters $a\,$.  In particular, $\Phi$ takes the zero locus
$F_2^{-1}(0)$ to $F_1^{-1}(0)$, and $\psi$ takes the parameter-space
projection of $F_2^{-1}(0)$ to that of $F_1^{-1}(0)$.
\medskip

In the case when $A(x)$ or $A(a,x)$ is the identity matrix (hence can be deleted
from the expressions~(\ref{e:equiv}) or~(\ref{e:Kequiv})) the maps or
deformations are {\em  right equivalent} (${\mathcal R}$-{\em
  equivalent}).  It is easy to check that if two real-valued functions are 
${\mathcal R}$-equivalent (up to a constant) then their gradient
vector fields are $\K$-equivalent.
\medskip

Much of the power of singularity theory resides in the central notion
of \emph{versal deformation}: a typical map germ has an
identifiable class of deformations from which \emph{all} others can be
obtained by suitable coordinate substitution.
\begin{defn}  \label{d:kversal}
The $r$-parameter deformation $F$ of $f$ is 
$\K$-{\em versal} if, given any $q$-parameter deformation $G(c,x)$ of
$f\,$, there is a map germ $\alpha:\br^q,0\to\br^r,0$ with the
property  that $G$ is $\K$-equivalent to the deformation
\[
\alpha^*F:(c,x)\mapsto F(\alpha(c),x).
\]
The deformation $F$ is \emph{$\K$-miniversal} if there exist
no $\K$-versal deformations of $f$ with fewer parameters.
\end{defn}
In other words, the parameters $a_i$ in the $\K$-versal deformation $F$ can
be expressed as suitable (smooth) functions $\alpha_i$ of the given 
parameters $c_j$ in order to recover $G$ from $F$ up to $\K$-equivalence.
Note that in general $\alpha$ is not a diffeomorphism, since (apart
from other considerations) $q$ and $r$ need not be equal.
\medskip

\rem It is a consequence of the general theory that if $F$ is
$\K$-miniversal then $F(a,\cdot)$ and $f$ are 
$\K$-equivalent germs if and only if $a=0\in\br^r$.
\medskip

If $f$ is nonsingular (that is $Df(0)$ has maximal rank) then the IFT
implies that every deformation of $f$ is $\K$-versal. In what follows
we therefore suppose that $f$ is singular.

Not every $f$ necessarily has a $\K$-versal deformation, but the failure
to do so is (in a precise sense) exceptional.  Moreover, and most
importantly, $\K$-versal deformations when they exist can be chosen to take
the explicit form
\begin{equation}  \label{e:vform}
F(a,x)\equiv a_1v_1(x)+\cdots+a_rv_r(x)+f(x)
\end{equation}
for a certain (nonunique) collection of polynomial map germs
$v_i:\br^n\to\br^p, 1\le i\le r$,
that can be constructed from $f$ by a standard algebraic procedure. \\

From Definitions~\ref{d:deform} and~\ref{d:kversal} we deduce the following:
\begin{prop}  \label{p:standard}
The zero set $Z_F:=F^{-1}(0)$ is a `standard object' in
$\br^r\times\br^n$ associated with $f\,$: 
the zero locus $Z_G$ of \emph{any} given deformation $G$
of $f$ as above is (up to diffeomorphism $\Phi$) the inverse image 
of $Z_F$ under the map $(c,x)\mapsto(\alpha(c),x)$.  
\endproof
\end{prop}
\medskip

We are concerned with bifurcations, so the points of $Z_G$ of
particular interest to us are those where the projection of $Z_G$ into
parameter space $\br^q$ fails to be a submersion: these are
points where the local structure of the set of solutions $x$ to
$G(c,x)=0$ can change as $c$ varies locally.
\begin{defn} \label{d:sing}
The \emph{singular set} $\Sigma_G$ of the deformation 
$G:\br^q\times\br^n\to\br^p$ is defined as
\[
\Sigma_G=\{(c,x)\in\br^q\times\br^n: \rank D_xG(c,x)<p\}.
\]
The \emph{discriminant} $\Delta_G$ is the projection of $Z_G\cap\Sigma_G$
into the parameter space $\br^q$.
\end{defn}
Thus $\Delta_G$ is the set of $c\in\br^q$ for which there exists some
$x\in\br^n$ with $G(c,x)=0$ but $D_xG(c,x)$ has rank less than $p$.
\medskip

The importance of $\Sigma_G$ for us lies in the following key result
much used in deformation (unfolding) theory.  See e.g.~\cite[p.49]{GI}
or~\cite[Ch.11,(18)]{CH} for a proof.
\begin{prop}
If $G$ is a submersion, so $Z_G$ is a smooth submanifold of
$\br^q\times\br^n$ of codimension $p\,$, then $\Sigma_G$ is the set of
singular points of the projection
$Z_G\subset\br^q\times\br^n\to\br^q$.
\endproof
\end{prop}

\rem 
If $n<p$ (a case important to us) then $\Sigma_G=Z_G$ and
the rank condition is vacuous.
\medskip

From Definition~\ref{d:kversal} it follows that if $F$ is a $\K$-miniversal
deformation of $f$ then, for any other deformation $G$ of $f\,$,
the discriminant $\Delta_G$ is (up to diffeomorphism $\psi$) the
inverse image of the discriminant $\Delta_F$ under the parameter map $\alpha$.
The discriminant $\Delta_F$ is also a \lq standard object'
and for polynomial $f$ it is given by a finite set of polynomial
equations in $\br^r$. 
\medskip

The policy statement below is the ideological basis for our
 geometrical treatment of bifurcation from $S$ in this paper:
\begin{quote}
\emph{The key to understanding the geometry of the zero locus of the
  deformation} $G$ \emph{of} $f$ \emph{is to understand the geometry of the
  associated map} $\alpha$ \emph{from the given parameter space of} $G$
  \emph{into that of a} $\K$-\emph{miniversal deformation} $F$ \emph{of} $f$, 
  \emph{and in particular the
  intersection of the image of $\alpha$ with the discriminant}
  $\Delta_F\,$.  \emph{Generic assumptions about this intersection lead to
  conclusions about generic geometry of} $\Delta_G$.
\end{quote}
Because of the degeneracy along $S$ we shall need to take extra
care in implementing this policy.
\subsection{Deformation along a manifold.}  \label{ss:defman}
There is a vital technical issue to deal with, before
we can exploit the above ideas. 
In our context as described in Section~\ref{ss:interpret}
the full parameter space is
$\br^q\times S$, with
the $\eps$-deformation terms vanishing on $\ots$. However, the notion of
versal deformation as described above applies only in a \nhd of a single point
in parameter space, whereas we would like to apply it in a 
consistent way over a whole \nhd of $\ots$ in $\br^q\times S$.  

In standard proofs~\cite{AVG,GI} of the existence of $\K$-versal
deformations the requisite coordinate changes are constructed by integrating
suitably chosen vector fields, the existence of these vector fields
being established by using the Preparation Theorem for smooth function germs.  
However, as pointed out in \cite{AVG} (see Section 4.5, Example 4),
there is a parametrised version of the Preparation Theorem 
which holds where the parameter belongs to a closed interval $[0,1]$
and where (as in our case of uniform normal degeneracy) the map germ
to be unfolded is independent of the parameter. The proof of this theorem uses (as 
expected) a smooth partition of unity.  The same proof applies equally well to a
parameter belonging to any smooth compact manifold (there are no
topological obstructions in view of the linear nature of the
constructions).  
This enables us to work with the `global' parameter
$(\eps,x)\in \br^q\times S$ for $|\eps|$ small.   We can therefore
state the following deformation result, customised for our purposes:
\begin{prop}  
\label{p:main}
Let $h:\br^k,0\to\br^{d+k},0$ be a smooth map germ which has a
$\K$-versal deformation $H:\br^r\times\br^k\to\br^{d+k}$.
Let $W$ be a \nhd of $0\times S_0$ in $\br^q\times S\times \br^k$ and let
$F:W\to\br^{d+k}$ be a smooth map
such that $F(0,x,y)=h(y)$ for all $x\in S$.  Then there
exists a \nhd $W'\subset W$ of $0\times S_0$ and a smooth map
\begin{equation}  \label{e:ourPhi}
\Phi:W'\to  \br^r\times\br^k :
                      (\eps,x,y)\mapsto (\psi(\eps,x),\phi(\eps,x,y)),
\end{equation}
as well as a nonsingular matrix $A(\eps,x,y)$ of size $d+k\,$ and
defined on $W'$, such that
\begin{enumerate}
\item
the identity
\begin{equation}
F(\eps,x,y) \equiv A(\eps,x,y)\,H\circ\Phi(\eps,x,y)
\end{equation}
holds for all $(\eps,x,y)\in W'$;
\item
the map $\br^k\to\br^k:y\mapsto\til y=\phi(\eps,x,y)$ is a diffeomorphism for
each $(\eps,x)\in \br^q\times S$, with $\phi(0,x,y)=y\,$ for all $x\in S\,$;

\item
$\psi(0,x)=0\in\br^r$ for all $x\in S\,$.

\end{enumerate}
\endproof
\end{prop}
Putting together Proposition \ref{p:main} and the
expression (\ref{e:vform}) we now obtain the deformation result we require
for application to the reduced bifurcation problem $\til F_\eps=0$ as in~(\ref{e:red}).  
\begin{cor} \label{c:global}
Let the vector field $F_0$ satisfy uniform normal degeneracy on $S$ as
in Section~\ref{ss:uniform}, so that in coordinates $(x,y)$ on the
\nhd\ $U$ of $S_0$ in $K=S\times\br^k$ we have 
$\til F_0(x,y)\cong h(y)$ where
$\til F_\eps$ is the reduced map as in (\ref{e:red}).
If $h$ has $\K$-versal deformation of the form 
\[
H(a,y)\equiv a_1v_1(y)+\cdots a_rv_r(y)+h(y)
\]
then, up to premultiplication by a nonsingular matrix depending smoothly
on $(\eps,x,y)$, we can express $\til F_{\eps}$ in the form
\begin{equation} 
\label{e:nform}
\til F_{\eps}(x,y) =  \til F(\eps,x,y) 
 \equiv a_1(\eps,x)v_1(\til y)+\cdots+a_r(\eps,x)v_r(\til y) + 
                                           h(\til y)\in\br^{d+k}
\end{equation}
where $y\mapsto\til y$ is a diffeomorphism germ at the origin in
$\br^k$ depending (smoothly) on $(\eps,x)\in \br^q\times S$ and where also 
\begin{equation} \label{e:azero}
a_i(0,x)\equiv 0
\end{equation}
for all $x\in S$ and $i=1,2,\ldots,r$.   \qed
\end{cor}
\subsection{Global and local geometry of the discriminant.} \label{ss:pullback}
\label{ss:pull-back}
We continue to assume uniform degeneracy, and 
without loss of generality we take $\til F_\eps$ to be of the
form~(\ref{e:nform}) where after the coordinate change $y\to\til y$ we
drop the tilde on $y\,$.
The identity (\ref{e:azero}) implies that we can write
\begin{equation}  \label{e:expand0}
a_i(\eps,x)=b_i(x)\eps + c_i(x)(\eps,\eps) +O(|\eps|^3) \in \br 
\end{equation} 
for $i=1,2,\ldots,r$, where $b_i(x)$, $c_i(x)$ are linear and bilinear maps 
$\br^q\to\br$ and $\br^q\times\br^q\to\br$ respectively.
Thus 
\[b_i(x)\,\eps=\sum_jb_{ij}(x)\,\eps_j\mbox{\quad and \quad}
c_i(x)(\eps,\eps)=\sum_{lm}c_{ilm}(x)\,\eps_l\,\eps_m\]
where 
\[
b_{ij}(x)=\dbd{a_i}{\eps_j}(\eps,x)|_{\eps=0}\mbox{\quad and\quad } 
c_{ilm}(x)=\frac{\partial^2 a_i}{\partial\eps_l\partial\eps_m}
                                                (\eps,x)|_{\eps=0}.
\]
We next write $\eps$ in `polar coordinates' as
\[
\eps=\rho s\,,\quad s\in \bs^{q-1},\quad \rho\in\br_+
\]
where $\br_+ =\{t\in\br:t\ge0\}$ and $\bs^{q-1}=\{s\in\br^q: |s|=1\}$.
From~(\ref{e:expand0}) the smooth map
\[
a:\br^q\times S\to\br^r:(\eps,x)\mapsto(a_1(\eps,x),\ldots,a_r(\eps,x))
\]
satisfies
\begin{equation}  \label{e:rho_eq}
a_i(\eps,x)=\rho b_i(x)s + \rho^2 c_i(x)(s,s) + O(\rho^3) \in \br
\end{equation}
for $i=1,\ldots,r$. We study the geometry of the map
\begin{equation}  \label{e:bdef}
b:\bs^{q-1}\times S\to\br^r:(s,x)\mapsto (b_1(x)s,\ldots,b_r(x)s).
\end{equation}
From this we obtain the geometry of the map
$a:\br^q\times S\to\br^r$
to first order in $\rho$ as the \emph{cone} on that of $b\,$.
We anticipate that under suitable generic assumptions on 
$b$ this cone structure will enable
us to capture the true geometry of the discriminant 
\[
\Delta_{\til F}:=a^{-1}(\Delta_H)\in \br^q\times S
\]
of $\til F$ near the origin. 

To investigate this closely we look at rays of the cone and their
intersections with 
$\Delta_H$ as $\rho\to 0$. It will be convenient to do this using a
`blow-up' or (in reverse) pinching technique as follows.  
First define the {\em pinching map}
\[
\tau:\br\times \bs^{q-1}\to \br^q : (\rho,s)\mapsto(\rho s)
\]
and write
\[
\til\tau:=\tau\times{id}:\br\times \bs^{q-1}\times S\to \br^q\times S\,,
\]
noting that since $\tau(\rho,-s)=\tau(-\rho,s)$ we may restrict
attention to $\rho\in\br_+:=\{r\in\br:r\ge0\}$ if convenient.
Next, let
\[
D_{\til F}^1:=b^{-1}(\Delta_H)\subset \bs^{q-1}\times S.
\]
Generic assumptions on $b$ will enable us to characterise the
local structure of $D_{\til F}^1$ which then determines
local models for the structure of
\[
\Delta_{\til F}^1:=
      \alpha^{-1}(\Delta_H)\subset\br_+\times\bs^{q-1}\times S
\]
where
\[
\alpha:= a\circ\til\tau:\br_+\times \bs^{q-1}\times S \to \br^r.
\]
(Note that here $\alpha$ plays the role of the map $\alpha$ in 
Definition~\ref{d:kversal}.)
Finally, the application of $\til\tau$ gives the local geometry of 
$\Delta_{\til F}=a^{-1}(\Delta_H)$.
\medskip

The following elementary fact about $\tau$ will be important.
\begin{prop}  \label{p:pinch}
Let $\gamma:\br,0\to\br\times\bs^{q-1},(0,x_0)$ be (the germ of) a smooth  
path such that
\[
\gamma(t)=(pt^m+O(t^{m+1}),s_0+tu+O(t^2))
\]
where $p\ne0\in\br$ and $u\ne0\in T_{s_0}\bs^{q-1}$ (so $u\perp s_0$
in $\br^q$).  Then the image of $\delta:=\tau\circ\gamma$ in $\br^q$ is (close
to the origin) the union of two $C^1$ arcs from the origin, each
having contact of order $(m+1)/m$ with their common tangent direction $s_0\,$.
\end{prop}
\proof  
Immediate, since $\tau\gamma(t)
=p\,t^m\,s_0+O(t^{m+1})s_0+ p\,t^{m+1}\,u + O(t^{m+2})$, 
so taking the parametrisation $t'=t^m$ gives
\[
\tau\gamma(t')-p\,s_0\,t'=p\,{t'}^{(m+1)/m}u 
\]
to lowest order in $t'$.
\endproof
\medskip

Observe that if $m$ is even then $\delta$ has a cusp at the origin in $\br^q$,
while if $m$ is odd then $\delta$ is a $C^1$ $1$-manifold passing through
the origin.
\medskip

The following geometrical construction will also play a key part in our investigations. 
\medskip

Let $H:\br^r\times\br^n\to\br^p$ be a $\K$-miniversal
deformation of the singular germ $h:\br^n,0\to\br^p,0$.
\begin{defn}
A vector $v\in\br^r$ is called a \emph{positive tangent vector} to
$\Delta_H$ at the origin if there exists $t_0>0$ and a differentiable
path $\gamma:[0,t_0)\to\br^r$ with $\gamma(0)=0\in\br^r$ 
such that
\begin{enumerate}
\item[(i)] $\gamma(t)\in\Delta_H$ for all $t\in[0,t_0)$
\item[(ii)] $\gamma\,'(0)=v$.
\end{enumerate}
The \emph{positive tangent cone} $T_H$ for
$\Delta_H$ at the origin is the union of all the positive 
tangent vectors at the origin.
\end{defn}
We now turn to applying this machinery to the reduced bifurcation
problem~(\ref{e:red}).
\subsection{Branch points and bifurcation}
As $\eps\in\br^q$ moves away from the origin in a given direction we
expect branches of solutions to $\til F_\eps=0$ to emanate from
certain points $(x,0)\in S_0\,$.
\begin{defn}
A \emph{solution branch} for $\til F$ at $x\in S$ is a differentiable
path $z:[0,t_0)\to Z_{\til F}\subset \br^q\times S\times\br^k$
for some $t_0>0$ such that 
\begin{itemize}
\item[(i)]
$z'(0)\ne0$ 
\item[(ii)]
$z(0)=(0,x,0)\in\br^q\times S_0$
\item[(iii)]
$z(t)\in\br^q\times S_0$ only if $t=0$.
\end{itemize}
A solution branch $z$ has \emph{direction} $s\in\bs^{q-1}$ if the
projection of $z'(0)$ into $\br^q$ is nonzero and in the direction
$s\,$; otherwise we say the branch $z$ is \emph{tangent} to $S$.
If a solution branch exists at $x\in S$ with direction $s$ then we say
that $x$ is a {\em branch point} with {\em branch direction} $s\,$. 
\end{defn} 
Recall (see the Remark after Definition~\ref{d:sing}) that since $k<d+k$ the
zero set $Z_{\til F}$ is here the same as the 
singular set $\Sigma_{\til F}$ and projects to 
$\Delta_{\til F}=a^{-1}(\Delta_H)\subset\br^q$.
\msk

From the geometry of $\Delta_H$ we obtain a necessary condition for branching.
\begin{prop}   \label{p:solns}
A necessary condition for $x_0\in S$ to be a branch point with branch
direction  $s_0\in \bs^{q-1}$ is that $b(s_0,x_0)$ lie in the
positive tangent cone $T_H\,$
for the discriminant $\Delta_H$ of the $\K$-versal deformation $H$ of $h\,$.
\end{prop}
\proof
Let $z$ be a solution branch at $x_0\in S$ with direction $s_0\,$
so that 
\[
z(t)=(\til z(t),y(t)) = (\rho(t)s(t),x(t),y(t)) \in \br^q\times S\times\br^k
\]
where $\rho(0)=0\,$, $\rho'(0)\ne0$ and $s(0)=s_0\,,x(0)=x_0$.
From~(\ref{e:rho_eq}) we have
\begin{equation}
\frac{d}{dt}a(\til z(t))\big|_{t=0}=\frac{d}{dt}\rho(t)b(s(t),x(t))
+O(\rho(t)^2)\,\big|_{t=0} = \rho'(0)b(s_0,x_0)
\end{equation}
and so $b(s_0,x_0)\in T_H$ since $a(\til z(t))\in\Delta_H$ for $t\in[0,t_0)$.  
\endproof 
\bigskip

In the examples considered below
the set $T_H$ turns out to be either a
hyperplane or a half-hyperplane.  Generically we expect the map
$b:\bs^{q-1}\times S\to\br^r$ to be transverse to $T_H\,$ 
in an appropriate sense, and the
transversality allows us to establish the existence of branch points
and to describe the local solution structure regardless of the 
terms $O(\eps^2)$ in $a(\eps,x)$.
\medskip
 
Because of the technical complications of the relevant singularity
theory, the approach described above is useful in practice
only for small values of $d$ and/or $k\,$. However, the case
$k=1$ is most significant, as it represents the least departure from 
normal nondegeneracy of $S\,$, while the case $d=1$ naturally occurs 
in the study of periodic orbits of
differential equations where there is an automatic circle group
action on the space of periodic functions (see Section~\ref{s:apps}).
Therefore the simplest but most important case $d=k=1$ merits detailed 
investigation.
\section{The fundamental case $d=k=1$.}  \label{s:d=k=1}
Here we have a $1$-dimensional equilibrium manifold $S$ of 
constant corank $1$. Throughout this section 
we continue to assume that the normal degeneracy is uniform;
some branching analysis applying to the general case is given below in
Section~\ref{s:q=k=1}. 
\msk

To simplify notation we drop the tilde throughout and write~(\ref{e:nform})
as
\begin{equation} \label{e:dk1}
   F_\eps(x,y)=F(\eps,x,y)\equiv g(\eps,x,y) + h(y) \in \br^2
\end{equation}
with multiparameter $\eps\in\br^q$ and $u=(x,y)\in S\times\br\,$ since
here $k=1$ and $d+k=2$.  
Each component of $h(y)$ is without constant or linear terms because of the
normal degeneracy of the zero set $S_0$.  
Supposing that at least one component has a
nonzero Taylor series at $0$, it is a straightforward exercise (see
\cite{GI} for example) to check that
$h$ is $\K$-equivalent to the germ $y\mapsto(0,y^m)$ for some integer $m>0$.  
We first take the least degenerate case $m=2$ and
study this example in detail in the geometric framework outlined above.
We then show how this approach extends to general $m$ and illustrate
in particular how it recovers and sets in a geometric context the
analysis of this problem for $q=2$ given by Hale and Taboas~\cite{HT2}.   
\subsection{The case $m=2$ \textrm{:} quadratic degeneracy.}  \label{ss:m=2}
Here
\[
F(\eps,x,y) = g(\eps,x,y) + h(y) = g(\eps,x,y) + (0,y^2) \in\br^2 
\]
for $(x,y)\in S\times\br$ with $\,\dim S=1\,$, and
$\eps\in\br^q$.
\subsubsection{Geometry of the zero set.}  \label{sss:geom}
\begin{figure}[htbp] 
   \centering
   \includegraphics[scale=0.3]{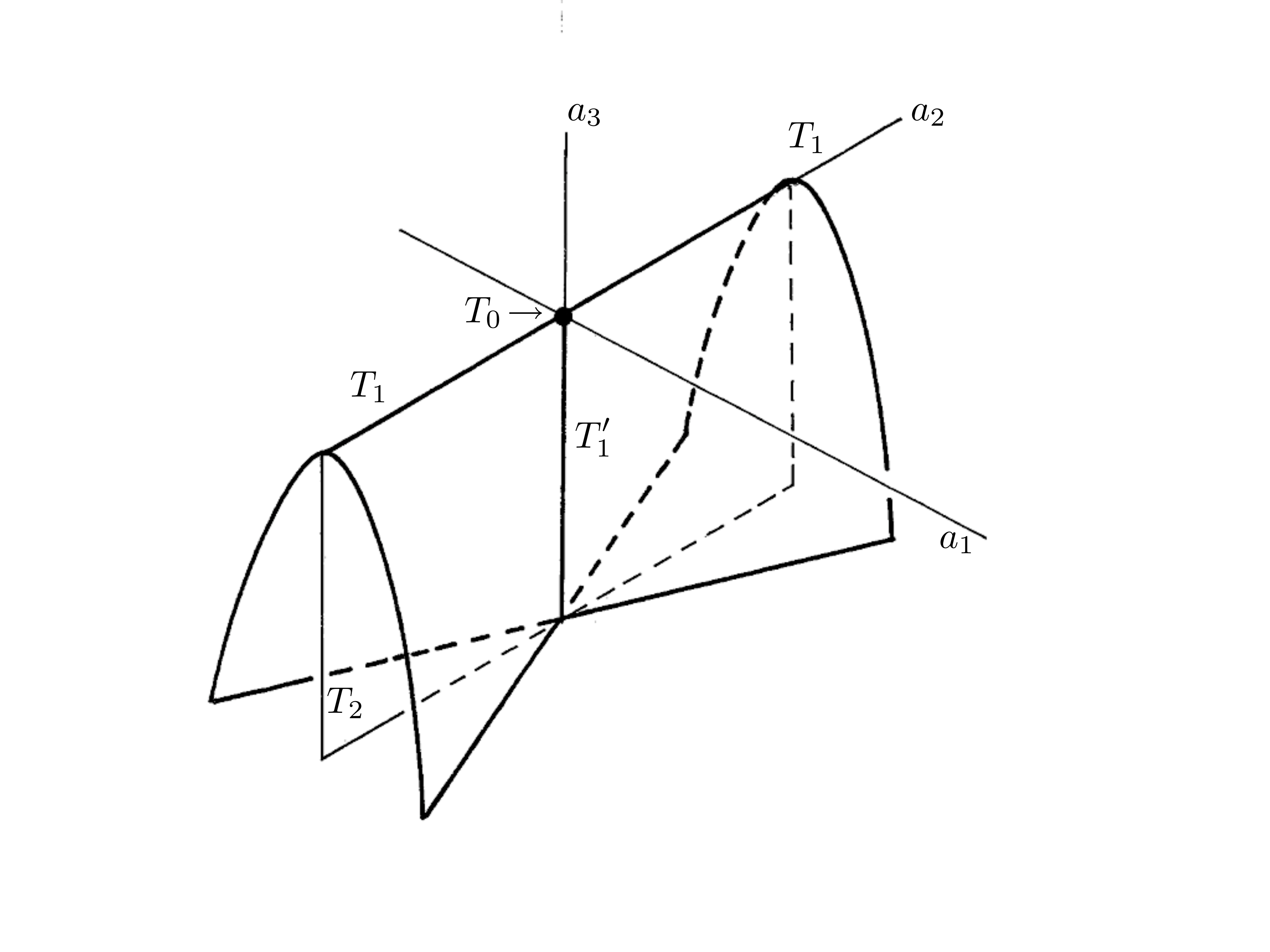} 
   \caption{Whitney umbrella: locus $a_1^2+a_2^2\,a_3=0$ showing the 
   stratification $\{T_0,T_1,T'_1,T_2\}$.}
   \label{fig:BB}
\end{figure}

A $\K$-miniversal deformation of $h(y)\equiv(0,y^2)$ is
$$H(a,y)=H(a_1,a_2,a_3,y)=(a_1+a_2y,a_3+y^2)$$ 
(see~\cite{AVG},\cite{GI}) and so by Corollary~\ref{c:global}
up to $\K$-equivalence we can suppose $F$ to have the form
\begin{equation} \label{e:unfy2}
F(\eps,x,y)=(a_1(\eps,x)+a_2(\eps,x)y\,,a_3(\eps,x)+y^2)
\end{equation}
for $|y|,|\eps|$ sufficiently small and for all $x\in S$, 
where $a(0,x)=0\in\br^3$ for all $x\in S$.  
As in Section~\ref{ss:pullback}, in order to
describe the bifurcation behaviour of $F$ we need to understand
the geometry of the discriminant
$\Delta_F=a^{-1}(\Delta_H)\subset\br^q\times S$
where (as is easily seen by eliminating $y$ from $H(a,y)=0$) the
discriminant $\Delta_H$ is given explicitly by
\begin{equation}  \label{e:delH}
a_1^2+a_2^2a_3=0.
\end{equation}
This is the equation of a \emph{Whitney umbrella} in $\br^3$: 
see Figure~\ref{fig:BB}.
%
\medskip

From Section~\ref{ss:pullback}, the key intermediate step is to
study the geometry of
\[
D_F^1:=b^{-1}(\Delta_H)\in \bs^{q-1}\times S
\]
where $b:\bs^{q-1}\times S \to\br^3$ is as in~(\ref{e:bdef}).
Following Proposition~\ref{p:solns} we first
identify the set $b^{-1}(T)$ where $T=T_H$ is the positive 
tangent cone to $\Delta_H$ at the origin.  
In the present case 
\[
T=\{a\in\br^3:a_1=0,a_3 \le 0\}
\]
since any smooth path in $\Delta_H$ from the origin can be checked to
have contact of order at least $\frac32$ with the plane $a_1=0$.
\medskip

We now make the following generic hypothesis, which is at the heart of our
whole approach to bifurcating from a normally degenerate manifold:
\bigskip

\noindent(H1)\quad The map $b:\bs^{q-1}\times S \to\br^3$ is transverse to $T$.
\bigskip

The notion of transversality to a variety such as $T$ rather than to a 
smooth manifold assumes a \emph{stratification} of the variety, that is a 
suitable decomposition of the variety into finitely many
manifolds of varying dimension, the smaller comprising the boundary of
the larger.  See~\cite[Ch.3,\S1]{AGLV} for example.
In this case there is an
obvious and natural stratification $\{T_0,T_1,T_2\}$ where $T_0$ is the 
origin, $T_1$ is the $a_2$-axis with $a_2\ne0\,$,
 and $T_2$ is the part of the $(a_2,a_3)$-plane
with $a_3<0\,$: here $\dim T_i=i$ for $i=0,1,2$.  
However, we see from
Figure \ref{fig:BB} that since our ultimate interest is in geometry 
of $\Delta_F$ we
should also take account of the fact that the negative $a_3$-axis $T_1'$ is a
natural $1$-dimensional stratum in $\Delta_H$.  We hence 
include this half-line $T_1'$ as a $1$-dimensional stratum for $T$
when interpreting the hypothesis (H1).  
%
\begin{figure}[htbp] 
   \centering
   \includegraphics[scale=0.5]{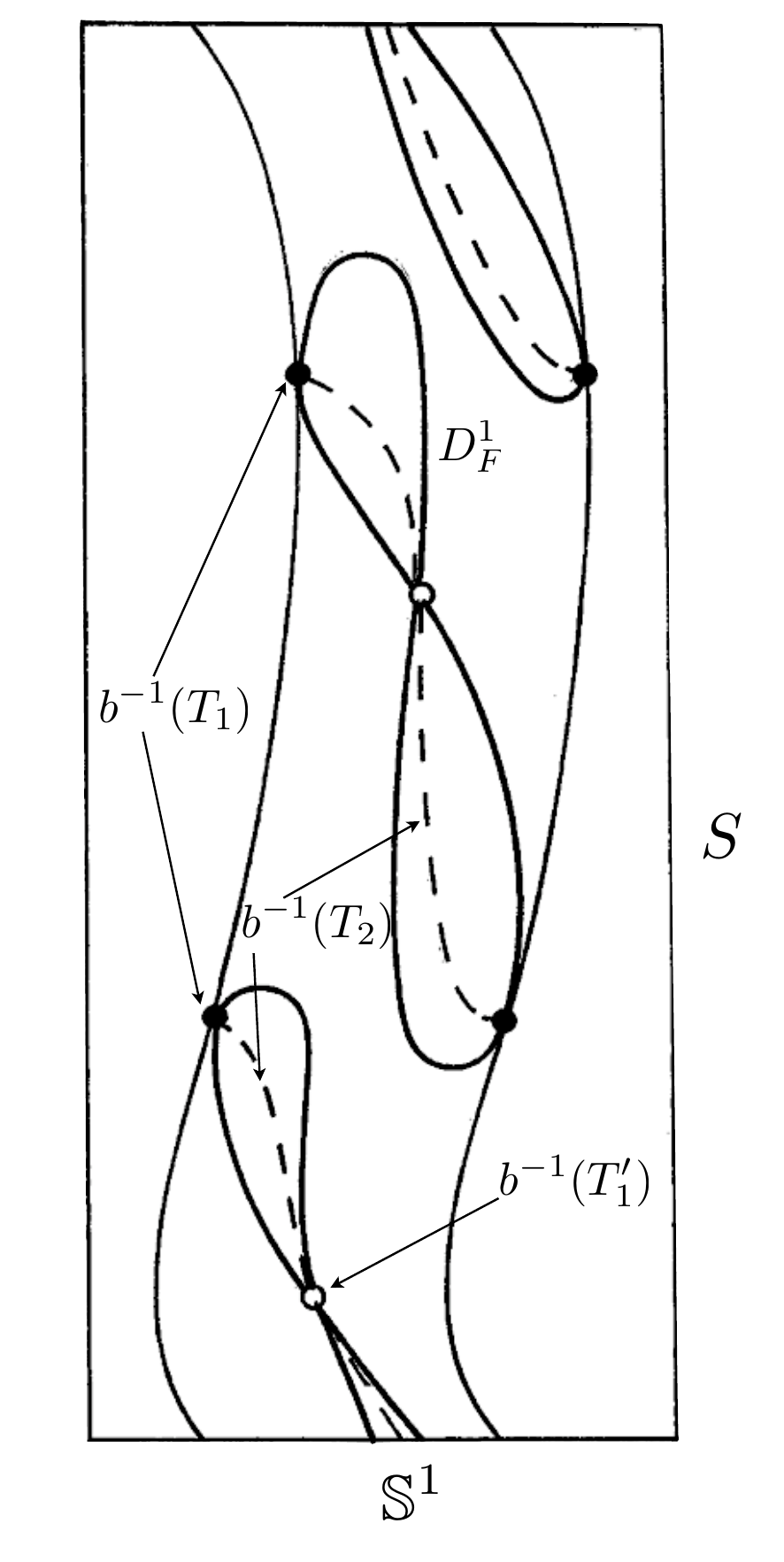} 
   \caption{Schematic representation of $b^{-1}(T)$ (dashed curve) and 
   $D_F^1=b^{-1}(\Delta_H)$ (solid curve) in the case $d=k=1$, $n=q=2$ (Whitney
     umbrella) showing end points (black dots) and intersection
   points (white dots). The whole of $D_F^1$ lies in a strip of width
   $\pi$ between 
   two parallel curves and has quadratic tangency with those at the
   end points. The $1$-dimensional manifolds $S$ and $\bs^1$ are each
   represented as closed intervals.}
   \label{fig:CC}
\end{figure}
\medskip

Since $\Delta_H$ and $T$ are $2$-dimensional objects in $\br^3$ the essential 
geometry of transversality is captured by taking $q=1$ or $2$ so that 
$\dim (\bs^{q-1}\times S) =1$ or $2$.  We first study carefully the case
$q=2$, and comment on the cases $q=1$ and $q\ge3$ later. 
\subsubsection{The `typical' case $q=2$.}  \label{sss:q=2}
\quad We take $S$ compact for ease of exposition; adjustments to the
noncompact case are straightforward.

The transversality hypothesis (H1) implies that
$b^{-1}(T)$ consists of a smooth $1$-manifold with boundary, thus
a finite collection of closed arcs or embedded circles
in $\bs^1\times S$.  The set of \emph{end points} of these 
arcs is $b^{-1}(T_1)$, while \emph{interior points} constitute
$b^{-1}(T_2)$ with possible {\em intersection points} $b^{-1}(T_1')$.
Note from the geometry of $\Delta_H$ and linearity of $b(s,x)$ in $s$
that if $(s,x)$ is an end point 
then so is $(-s,x)$.  Moreover, the intersection of
$S^1\times\{x\}$ with $\Delta_F^1$ lies in the semicircular arc
$\{s\in S^1:b_3(s,x)\le0\}$ since all of $\Delta_H$ lies in the
half-space $a_3\le0$.
A schematic representation of this geometry is shown in Figure \ref{fig:CC}. 
%
\begin{figure}[htbp] 
   \centering
   \includegraphics[scale=0.5]{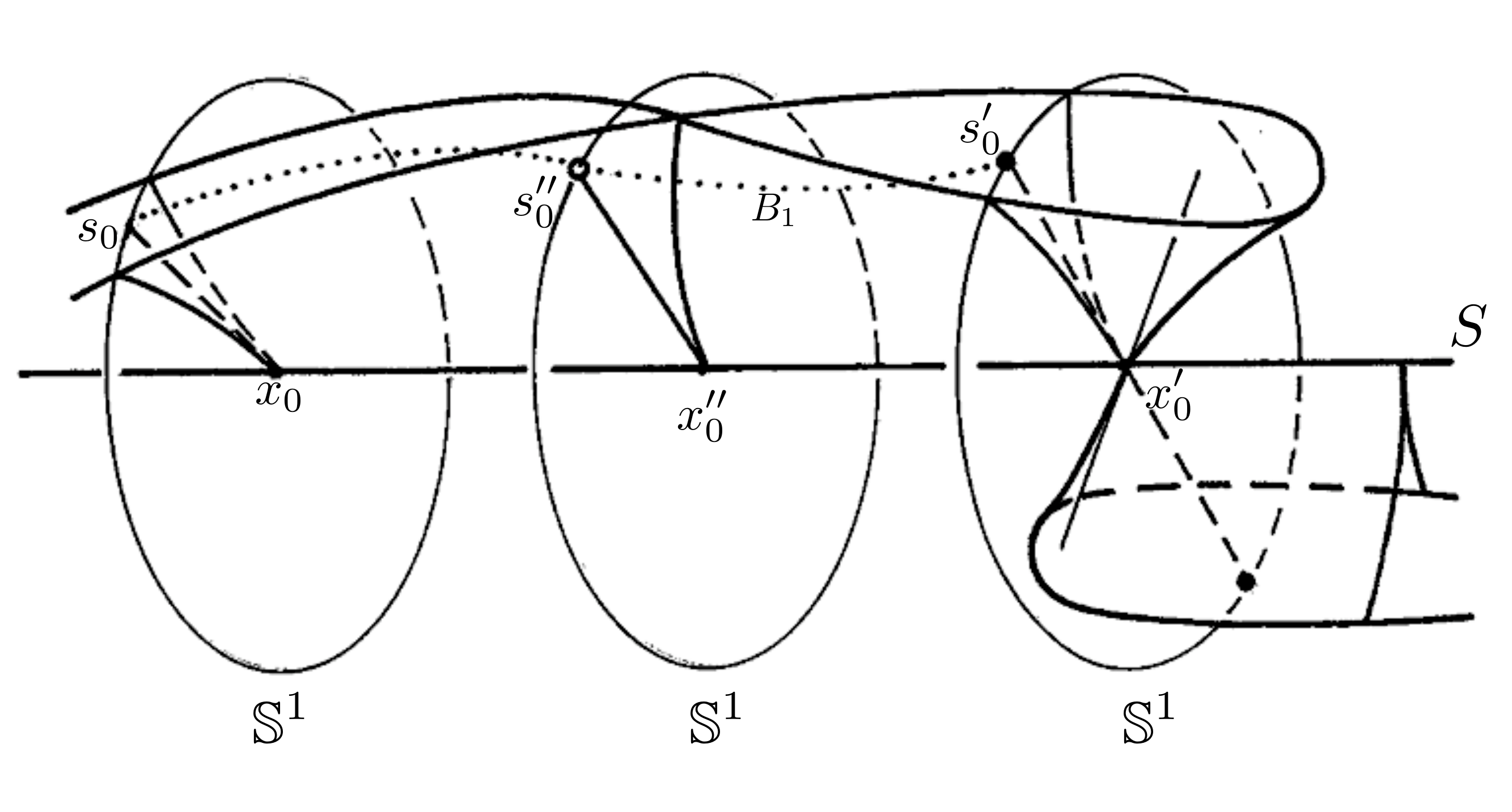} 
   \caption{Schematic diagram of the discriminant $\Delta_F$ with
     cusped ridge along $0\times S$, showing the structure close to
     $x_0,x_0',x_0''$ on $S$ determined by an interior point
     $(s_0,x_0)$, end point $(s_0',x_0')$ and intersection point
     $(s_0'',x_0'')$ (indicated simply by $s_0$ etc.) respectively.}
   \label{fig:DD}
\end{figure}
\medskip

For $(s,x)$ in a local coordinate chart around
$(s_0,x_0)\in\bs^1\times S$ write  $(s,x)=(s_0+v,x_0+w)$ 
where $(v,w)\in \br\times\br$.  Using~(\ref{e:rho_eq}) we have
\begin{equation} \label{e:alpha}
\alpha(\rho,v,w)=a\circ\til\tau(\rho,s,x)
                =\rho b(v,w)+\rho^2 c(v,w)+O(\rho^3)\in\br^3
\end{equation}
where now $b(v,w)$ is shorthand for $b(s,x)$ and analogously
$c(v,w)$ stands for $c(x)(s,s)$.  
From~(\ref{e:delH}) we see $\alpha(\rho,v,w)\in\Delta_H$ when
\begin{equation} \label{e:subindelta}
\bigl(\rho b_1+\rho^2c_1\bigr)^2
	+\bigl(\rho b_2+\rho^2c_2\bigr)^2\bigl(\rho b_3+\rho^2c_3\bigr)
				+O(\rho^4)=0
\end{equation}
which after cancelling $\rho^2$ reduces to
\begin{equation}  \label{e:psi}   
\beta(\rho,v,w):=b_1^2+(2b_1c_1+b_2^2b_3)\rho+d(\rho,v,w)\rho^2 = 0
\end{equation}
where $b_i,c_i$ are components of $b(v,w),c(v,w)$ and where 
$d(\rho,v,w)$ is a smooth real-valued function.  Now we study three 
possibilities, illustrated in Figure \ref{fig:DD}. 
\bigskip

(\emph{i})\quad $(s_0,x_0)$ is an interior point: $b_1=0$ and
 $b_2b_3\ne0$.
\medskip

The IFT at $(v,w)=(0,0)$ implies that if 
$\til\rho:=b_1^2-\beta$
then the map 
$(\rho,v,w)\mapsto(\til\rho,v,w)$ is a local diffeomorphism at the origin. In 
$(\til\rho,v,w)$ coordinates the locus $\beta=0$ is given by 
$\til \rho=b_1^2$, which (since $(0,0)$ is a regular point of the
function $b_1$ by hypothesis (H1) of transversality to $T_2$) 
is a smooth surface in $\br^3$ quadratically tangent to the curve 
$b_1(v,w)=0$.

Finally, to reconstruct the geometry in the original $(\eps,x)$
deformation coordinates, we apply the pinching map $\til\tau$ 
which takes the surface $\beta=0$ locally to a surface with a
cusped ridge along $0\times S\,$; the tangent direction of the cusp at
$(0,x_0)$ is $s_0\in \bs^1$, while as $x$ moves away from
$x_0$ in $S$ the tangent direction $s$ of the cusp is given by the smooth curve 
\[
B_1=\{(x,s)\in \bs^1\times S: b_1(x,s)=0\}
\]
through $(s_0,x_0)$ in $\bs^1\times S$.  See Figure \ref{fig:DD}. 
%
\begin{figure}[htbp] 
   \centering
   \includegraphics[scale=0.7]{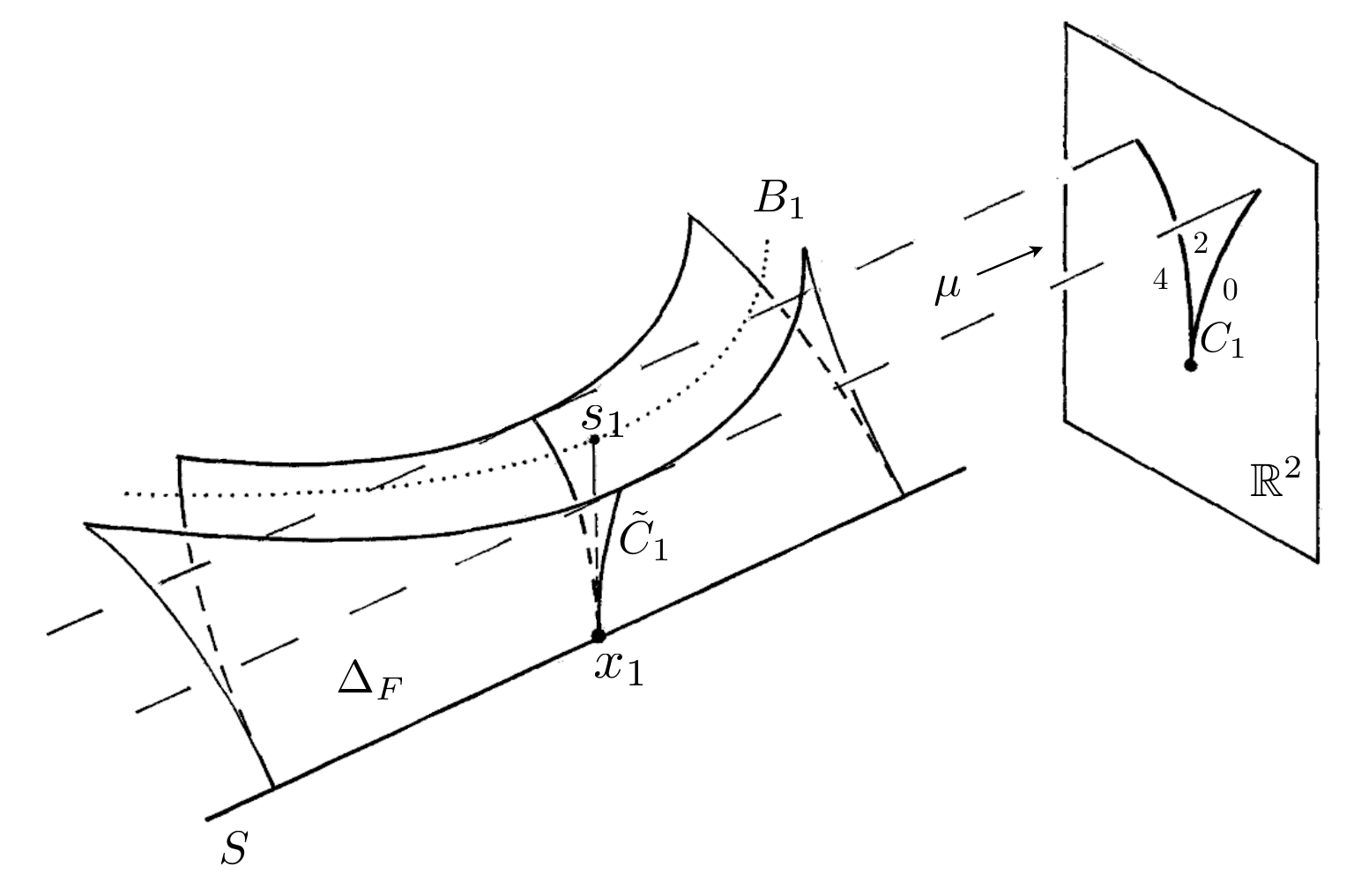} 
   \caption{The projection $\mu:\Delta_F\to\br^2$ close to
   $(0,x_1)\in\br^2\times S$  where $(s_1,x_1)$ 
   is an interior fold point of the curve $B_1$. 
   The cusp curve $\tilde{C}_1$ projects to 
   a cusp curve $C_1$ (the bifurcation set) in $\br^2$ tangent to the
   $s_1$ direction.  The numbers $0,2,4$ indicate the number 
   of solutions $(x,y)\in K$ to $F(\eps,x,y)=0$ close to $(x_1,0)\in K$ 
   for $\eps$ close to the origin in $\br^2$.}
   \label{fig:HH}
\end{figure}
\bigskip

(\emph{ii})\quad $(s_0,x_0)$ is an end point (denoted $(s_0',x_0')$): 
$b_1=b_3=0, b_2\ne0$.
\medskip

From the hypothesis (H1) of transversality to $T_1$ we have
\(
\dbd{(b_1,b_3)}{(v,w)}\ne0
\)
at $(0,0)$, and so the map
\[
\br^2\to\br^2 : (v,w)\mapsto (\til v,\til w):=(b_1,2c_1b_1+b_2^2b_3)
\]
is a diffeomorphism germ at the origin.  In new coordinates the function
$\beta$ takes the form 
$\beta(\rho,\til v,\til w)=Q(\rho,\til v,\til w) + O(\rho^3)$ 
where
\begin{equation}  \label{e:qf}
Q(\rho,\til v,\til w) =\til v^2 + \til w\rho + e\rho^2 
\end{equation}
which is a nondegenerate quadratic form regardless of the value
of the coefficient $e\,$. The locus $Q=0$ is a cone $C$ tangent to the
$(\til v,\til w)$-plane along the $\til w$-axis.  By the Morse Lemma
there is a local diffeomorphism (with derivative the identity) that
takes the locus $\beta^{-1}(0)$ to the cone $C$ and, moreover, 
(\ref{e:qf}) shows that $\beta^{-1}(0)$ is quadratically tangent to the
$(\til v,\til w)$-plane along the smooth curve $\til v=b_1=0$.  
Finally, by Proposition~\ref{p:pinch} the pinching map $\til\tau$ takes
$\beta^{-1}(0)$ to a surface with $\frac32$-power cusped ridge along 
$0\times S$ and terminating at $(0,x_0')$, reminiscent of a ship's prow: see 
Figure~\ref{fig:DD}.
%
\begin{figure}[htbp] 
   \centering
   \includegraphics[scale=0.5]{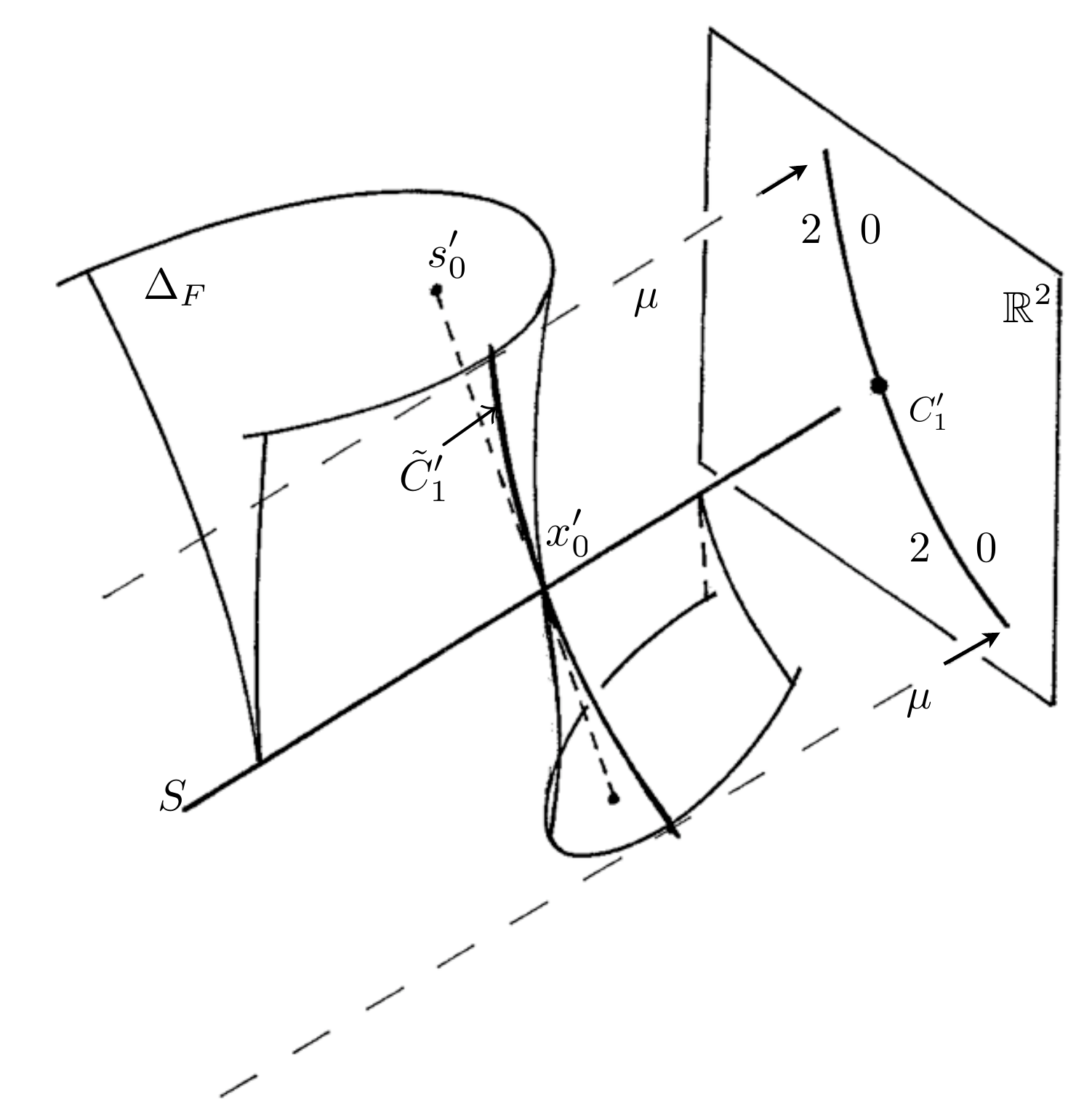} 
   \caption{The projection $\mu:\Delta_F\to\br^2$ close to
     $(0,x_0')\in\br^2\times S$ where $(s_0',x_0')$ 
   is an end point of the curve $B_1$. 
   The curve $\tilde{C}'_1$ projects to a pair of arcs $C'_1$ 
   (the bifurcation set) in $\br^2$ tangent to the $s_0'$ direction. 
     The numbers $0,2$ indicate the number of solutions  $(x,y)\in K$ to
     $F(\eps,x,y)=0$ close to $(x_0',0)\in K$ for $\eps$ close to 
     the origin in $\br^2$.}
   \label{fig:JJ}
\end{figure}
\bigskip

(\emph{iii})\quad $(s_0,x_0)$ is an intersection point (denoted $(s_0'',x_0'')$):
$b_1=b_2=0\,, b_3<0\,$.
\medskip

By the hypothesis (H1) of transversality to $T_1'$ the map
$(v,w)\mapsto(b_1,b_2)$ is a local diffeomorphism at $(0,0)$.
From~(\ref{e:delH}) with $a_3<0$ we can express $\Delta_H$ as
\[
(a_1+\ell a_2)(a_1-\ell a_2)=0
\]
where $\ell=\sqrt{|a_3|}$, which on cancelling $\rho^2$ becomes
$b_1\pm\ell b_2=0$ with $\ell=\sqrt{\rho|b_3|}$.
Hence for fixed $\rho>0$ the expression~(\ref{e:psi})
represents a pair of smooth curves crossing 
transversely at the origin, with tangent directions both tending to
that of $b_1(v,w)=0$ as $\rho\to0$.  See Figure \ref{fig:DD}.
\subsubsection{Bifurcation geometry when $q=2$.}  \label{sss:bifgeom}
We have described above the three generic types of local geometrical structure
found in the discriminant $\Delta_F\,$, where we recall that 
(since $1=d<p=d+k=2$) the discriminant is the projection into
$\br^2\times S$ of the zero set:
\[
\Delta_F=\{(\eps,x)\in\br^2\times S:\exists\, y\in\br\,, \,F(\eps,x,y)=0\}.
\]
However, it is important to keep sight of the fact that in terms of
the original problem~(\ref{e:red}),(\ref{e:dk1}) it is $(x,y)$ that
are the variables while $\eps$ is the perturbing parameter.  Therefore
to find the bifurcation behaviour of solutions $(x,y)$ as $\eps$ varies in
$\br^2$ we must also study the geometry of
the projection $\mu:\Delta_F\to\br^2$ into the $\eps$-plane and in
particular its singularities. The image of the
singular points (i.e. the set of singular \emph{values}) is the
\emph{bifurcation set}; this can also be interpreted as the \emph{apparent
outline} (or \emph{apparent contour}: see~\cite{BG},\cite{Dem}) of $\Delta_F$ 
when viewed along the $1$-manifold $S$.
\medskip

The analysis of the previous section showed that 
$\Delta_F$ has a cusped ridge along $0\times S$ with the direction of
the axis of the cusp (indicated on the unit circle $\bs^1\subset\br^2$)
given by the curve $B_1=b_1^{-1}(0)\subset \bs^1\times S$. 
(For any interior point  $x\in S$
there can be only one $s\in \bs^1$ with $(s,x)\in B_1$, so
cusps may impinge on $(0,x)\in0\times S$ from only one direction,
while at end points they do so from two opposite directions.)
Therefore the singular points of $\mu|B_1:B_1\to \bs^1$
play a key role in the structure of the bifurcation set for $F$.
Accordingly, we make the following generic hypothesis: 
\bigskip

\noindent(H2)\quad The map $\mu|B_1:B_1\to \bs^1$ has only {\em fold}
singularities, and these are at interior points of $B_1\,$.
\bigskip

A fold singularity of $\mu|B_1$ is characterised in local coordinates
$(v,w)$ on $\bs^1\times S$ by 
\begin{equation}  \label{e:bproj1}
\dbd{b_1}w=0,\quad \dbd{^2b_1}{w^2}\ne0 \quad \mbox{at}\quad (v,w)=(0,0).
\end{equation}
Following the classification in the previous section we now consider
in turn the three types of point on the curve $B_1$.
\bigskip

(\emph{i})\quad Interior points.
\medskip

The first result confirms what is evident from a simple sketch.  See
Figure~\ref{fig:HH}.
\begin{prop}  \label{p:bproj1}
Let $(s_1,x_1)$ be an interior point of $B_1\subset \bs^1\times S$ 
such that the projection $\mu|B_1:B_1\to\bs^1$ has a fold singularity at
$(s_1,x_1)$.  Then the apparent outline in $\br^2$ of a \nhd of
$(0,x_1)$ in $\Delta_F$ viewed along $S$ is a curve $C_1$ with a cusp 
of order $\frac32$ at the origin and with tangent direction $s_1\,$ there.  
The number of solutions $(x,y)$ to $F_\eps(x,y)=0$ increases by two as
$\eps$ crosses each arc of $C_1$ in the appropriate (and same) direction of
rotation about the origin. 
\end{prop}
\proof 
In local coordinates $(\rho,v,w)$ the apparent outline $\Gamma_1$ 
of $\Delta_F^1$ viewed in the direction of $S$ is  by definition given by
\[
\Gamma_1=\left\{(\rho,v)\in\br_+\times \br:\exists\, w\in S\,, \,\beta(\rho,v,w)
                                 =\dbd\beta{w}(\rho,v,w)=0\right\},
\]
and we have already seen in Section~\ref{sss:q=2} 
(\emph i) that in suitable coordinates (in which the role of $w$ is
unaffected) we may take
$\beta=b_1^2-\rho$. Hence $\dbd\beta w=2b_1\dbd{b_1}w$ which vanishes 
where $b_1=0$ (so then $\rho=0$ if $\beta=0$) and also where $\dbd{b_1}w=0$.  
By the IFT and the fold condition~(\ref{e:bproj1}), 
the latter locus is $\br_+\times B_1'$ where $B_1'$ is a
smooth curve through the origin in $\br^2$
not tangent to the $w$-axis.  The intersection $\Gamma_1$ of $\br_+\times B_1'$  
with $\Delta_F^1$ is a smooth curve quadratically tangent to the $(v,w)$ plane.  
By Proposition~\ref{p:pinch} the curve $\til C_1:=\til\tau(\Gamma_1)$ has a
$\frac32$-power cusp at $(0,x_0)$ and projects under $\mu$ to
a curve $C_1$ with the properties described. 
\endproof
\medskip

\emph{Note}: The geometry of the curve $\til C_1$ in
$\Delta_F=a^{-1}(\Delta_H)$ and the expression~(\ref{e:unfy2})
show that as $\eps$ crosses the arcs of $C_1$ 
one pair of solutions $(x,y)$ is created (or annihilated) close to
$y=y_1$ where $y_1=\sqrt{-a_3(\eps,x)}$, and the other pair close to $-y_1$.
\bigskip

(\emph{ii})\quad End points.  
\medskip

The analysis in Section~\ref{sss:q=2} showed that at an end-point 
the surface $\Delta_F^1=\beta^{-1}(0)$ has a Morse singularity, thus a local 
cone structure, with the cone tangent to the $(v,w)$-plane along the
curve $B_1$.  Let $\chi=\dbd\beta w\,$.
From~(\ref{e:qf}) we see that $\chi=0$ when $\rho=0$ and $b_1=0\,$, 
that is on the curve $B_1$. Differentiating the expression (\ref{e:psi})
with respect to $w$ we then observe that provided the generic condition
\bigskip

\noindent(H$2'$)\quad 
   $\dbd{}w (2b_1c_1+b_2^2b_3) \ne 0 $
\bigskip

\noindent holds at $(0,0)$ then $\dbd\chi\rho(0,0)\ne0$ and so by the IFT
the locus $\chi^{-1}(0)$ locally has the form of a graph of a smooth function
$\rho=\rho(v,w)$.  Now $\dbd{\rho}w=-\dbd{\chi}w/\dbd{\chi}\rho$ which
is nonzero at an end point because on $B_1$ we have
$\dbd{\chi}w=2\bigl(\dbd{b_1}w\bigr)^2\ne0$ by the hypothesis (H2).
Thus the smooth surface $\chi^{-1}(0)\subset \br\times\br^2$ 
cuts the $(v,w)$-plane and hence the
\lq cone' $\Delta_F^1$ transversely along $B_1$ (away from the origin)
and so must cut the cone again along another smooth curve $\Gamma_1'$
through the origin and transverse to the $v,w$-plane.
Then $\til C_1'=\til\tau(\Gamma_1')$ is a smooth curve whose image under
$\mu$ is by Proposition~\ref{p:pinch} a $C^1$ arc $C_1'$ through the
origin in $\br^2$.  See Figure~\ref{fig:JJ}. To summarise:
\begin{prop}
Let $(s_0',x_0')$ be an end point of $B_1\subset \bs^1\times S$ and
assume that (H$2'$) holds.  Then the apparent outline in $\br^2$ of a \nhd of
$(0,x_1')$ in $\Delta_F$ viewed along $S$ is a $C^1$ arc $C_1'$ through the
origin, across each branch of which the number of solutions
$(x,y)$ to $F_\eps(x,y)=0$ increases by two, taken in the appropriate 
(and opposite) directions of rotation about the origin. 
\endproof
\end{prop}

\rem The arc $C_1'$ is typically \emph{not} the same as the
`end-point' curve $a^{-1}(T_1)$ which itself corresponds only to
solutions with $y=0$.  Thus soon after a pair of solutions is created across
$C_1'$ one of them will pass through $y=0$, that is through $S_0\,$. 
However, we are here solving the reduced equation~(\ref{e:red}), so the
corresponding equilibria of~(\ref{e:orig}) for $\eps\ne0$ need not lie
on $S$ although they will have zero $K$-component and therefore lie in $L$.
\bigskip

(\emph{iii})\quad Intersection points.
\medskip

Intersection points are zeros of the map $(b_1,b_2):\bs^1\times
S\to\br^2$ and the zero set of $b_1$ is the curve $B_1\,$.  We make a
further generic hypothesis:
\bigskip

\noindent(H3)\quad The function $b_2$ does not vanish at the singular
points of $\mu|B_1\,$.
\bigskip

An intersection point $(s_0'',x_0'')$ corresponds to an arc of
self-intersection of $\Delta_F$ emanating from $(0,x_0)\,$: see Figure~\ref{fig:DD}.
This arc projects under $\mu$ to an arc $C_1''$ from the origin in
$\br^2$, and $(s_0'',x_0'')$ is a regular point of $\mu|B_1$ by the hypothesis (H3).
Thus no bifurcations occur across $C_1''$, but two solutions of the
reduced equation of opposite sign in $y$ pass through $y=0$ (that
is, through $S_0\,$) simultaneously as $\eps$ crosses $C_1''$, so the
corresponding solutions of~(\ref{e:orig}) pass through $L$.
\bigskip

We have now studied the three generic types of local
bifurcation behaviour, and we assemble them into a global statement.
Recall that $d=k=1$ (so the compact manifold $S$ is $1$-dimensional and is
degenerate with corank $1$) and the degree $m$ of degeneracy is $2\,$.
The perturbation parameter $\eps$ lies in $\br^2$.
\medskip

\begin{figure}[htbp]
   \centering
   \includegraphics[scale=0.7]{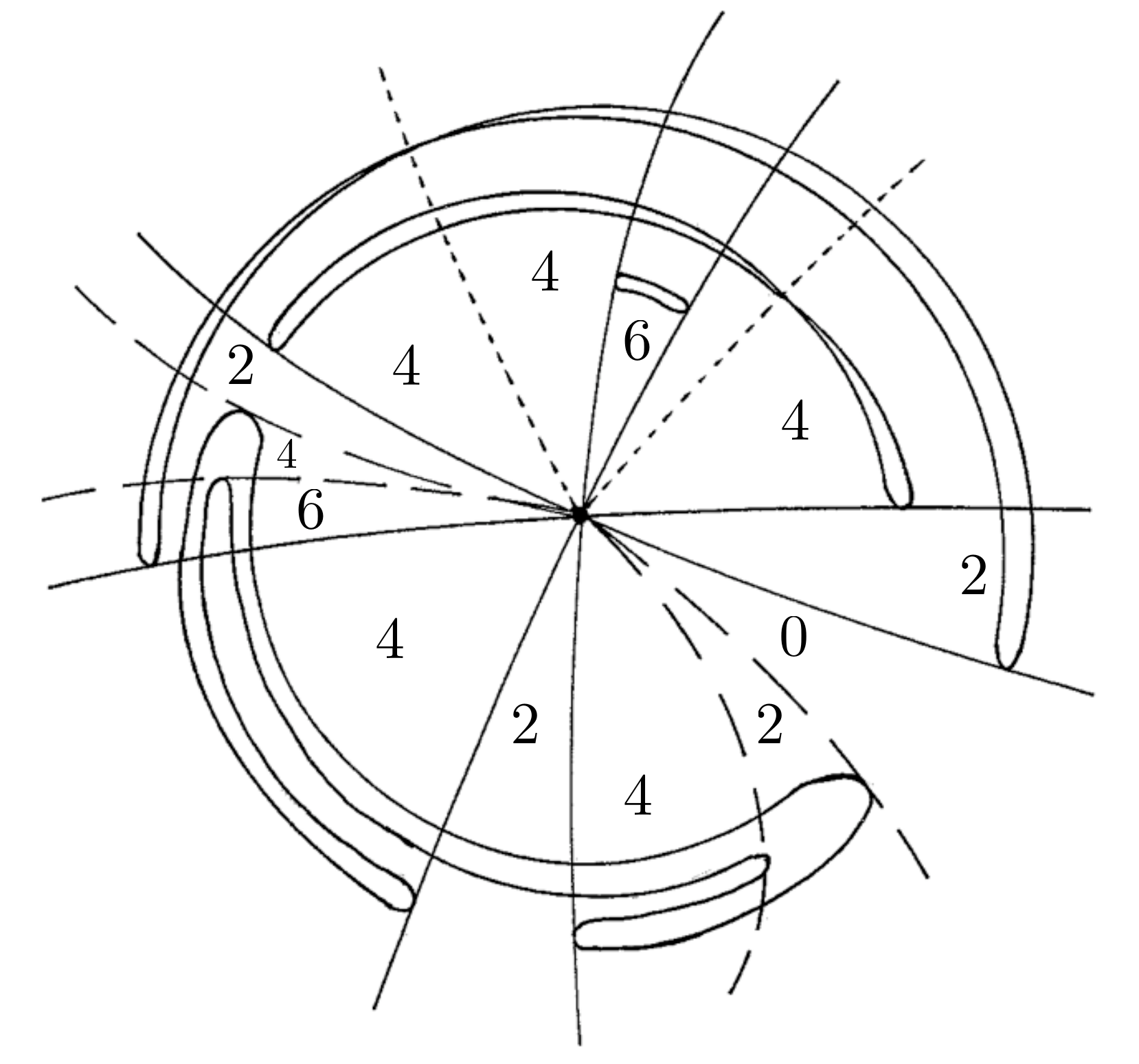} 
   \caption{Schematic bifurcation diagram  for $d=k=1,\, q=2$ and $m=2$
     (see text), showing end arcs (unbroken lines) and fold arcs (long
     dashes) through the origin in the parameter space $\br^2$. 
     Intersection arcs (short dashes) are also shown. Numbers
     indicate the number of solutions in each region. Note that for
     each opposite pair of fold arcs the two folds approach from the
     same side in the plane.}
   \label{fig:figbifA}
\end{figure}

Let $\kappa$ denote the local involution in a \nhd\ of $S$ in $M$ that
in terms of normal bundle coordinates $TS\oplus K\oplus L$ corresponds
to $(x,y,z)\mapsto(x,-y,z)$.
\begin{theo}  \label{t:bif1}
Under the generic hypotheses {\rm (H1),(H2),(H$2'$),(H3)} the bifurcation set
for~(\ref{e:orig}) consists of a finite set of $C^1$ arcs 
emanating from the origin in $\br^2$.  The arcs are of three types:
\emph{fold, end} and \emph{intersection} arcs. These correspond to the
following bifurcation behaviour as $\eps$ rotates about the origin in
$\br^2$ in (locally) the appropriate direction.
\medskip

\emph{Interior arcs}.\quad These occur in pairs with $\frac32$-power
cups at the origin.  As $\eps$ crosses the first arc a pair of solutions
is created away from $S$ by a saddle-node (fold) bifurcation, and
as $\eps$ crosses the second arc a second saddle-node pair of solutions is
created close to the $\kappa$-images of the first pair.
\medskip

\emph{End arcs}.\quad These occur in opposite pairs forming $C^1$
curves through the origin, and each arc is
accompanied by an arc of $a^{-1}(T_1)$ with the same tangent at the
origin.  As $\eps$ crosses an end arc a pair of solutions
is created away from $S$ in a saddle-node bifurcation; as $\eps$
crosses the companion arc one of these solutions passes through $L$.
This behaviour is repeated at the opposing end arc and with
the same sense of rotation about the origin.
\medskip

\emph{Intersection arcs}.\quad As $\eps$ crosses an intersection arc with
nonzero speed two solutions of the reduced equation
with $K$-components of opposite signs 
pass  through $L$ at the same point and with nonzero speed.
\end{theo}
This theorem in particular recovers and provides a geometric context
for the results of Hale and Taboas~\cite{HT2} for the case $m=2$.
It also extends those results by taking account of self-intersection 
points of the determining curve $B_1$. See Figure~\ref{fig:figbifA}, and
compare Figures~1 and~2 in~\cite{HT2}. 
%
\subsubsection{Cases $q\ne2$.}
For $q=1$ (that is $\eps\in\br$) the transversality hypothesis
(H1) implies that the image of $b:\bs^0\times S\to\br^3$ meets $T$ 
only at points of $T_2$ and that $0$ is a regular value of
$b_1:S\to\br$. Hence $b_1^{-1}(0)$ is a discrete set of points $\{(\pm1,x_i)\}$. 
(This may be visualised in Figure~\ref{fig:CC}
by taking a pair of vertical lines corresponding to fixed values  
$\pm s\in \bs^1$.)  To each $x_i$ with either $\eps>0$ or $\eps<0$  
there corresponds a pair of solution branches in
$\Delta_F\subset\br\times S$ emanating from $(0,x_i)$; 
the branches are such that locally $x-x_i$ varies as $|\eps|^{\frac32}$.
In particular this describes the branching behaviour for $q=2$ when $\eps$
leaves the origin along an arc {\em not} in a direction tangent to an
interior, end or intersection arc as in Theorem~\ref{t:bif1}.
\smallskip

For $q=3$ some interesting new geometry arises. For $x\in S$ the map
$a({\cdot},x):\br^3\to\br^3$ is typically nonsingular at the
origin, although generically there is a set of isolated points
$\{x_j\}\subset S$ where the derivative $D_j$ of $a({\cdot},x_j)$ has rank
$2$.  For $x\notin\{x_j\}$ and for sufficiently small
$|\eps|$ the discriminant $\Delta_F=a^{-1}\Delta_H$ 
intersects $\br^3\times\{x\}$ in diffeomorphic copy
$\{\Delta_H^x\}$ of $\Delta_H$ (more
precisely, a \nhd of the origin in $\Delta_H\,$, which amounts to the
same thing).  To capture the structure of the bifurcation set
we study the map $b:\bs^2\times S\to\br$ as in~(\ref{e:expand0})
and the geometry of $b^{-1}(T)$ where $T=T_H\,$.  Writing
$b_x(s)=b(s,x)$ we note that $B_x:=b_x^{-1}(T)=B_x^0\cup B_x^1$ 
where $B_x^0$ is a semicircular (open) arc $b_x^{-1}(T_2)$
of a great circle $C_x$ in $\bs^2$ with antipodal end-points
$B_x^1=b_x^{-1}(T_1)$, while $D_F^1(x):=b_x^{-1}(\Delta_H)$ is a figure
eight curve on $\bs^2$ passing through these end-points and with
self-intersection at $b_x^{-1}(T_1')$.  
As $\rho\to 0$ the intersection $\rho\bs^2\cap\Delta_H^x$  
flattens to $\rho b_x^{-1}(T)$ at rate $\rho^{\frac32}$.

In this setting the bifurcation set in $\br^3$,
which is the set of critical values of the projection
$\mu:\Delta_F\subset\bs^3\times S\to\br^3$, is conveniently
interpreted as (to first order in $\rho$) a cone on the 
{\em envelope} $E$ of the family of curves
$\bigl\{D_F^1(x)\bigr\}_{x\in S}$in $\bs^2$. The curve $E$ separates points
$s\in\bs^2$ according to the number $n_s$ of points $x\in S$ with $s\in
D_F^1(x)$.  The generic local geometry of envelopes of smooth curves
in $\br^2$ (or $\bs^2$) is well known (see~\cite{BG} for example): the
envelope consists of smooth fold curves across which $n_s$ changes by
$2$, meeting at isolated cusp points.  These local configurations are
stable, so that the same description applies to the envelope of the
family  $\bigl\{\rho\bs^2\cap\Delta_H^x\bigr\}_{x\in S}$ (where higher
order terms in $\rho$ are not neglected) for $\rho$ 
nonzero and sufficiently small.
%
\begin{figure}[htbp] 
   \centering
   \includegraphics[scale=0.7]{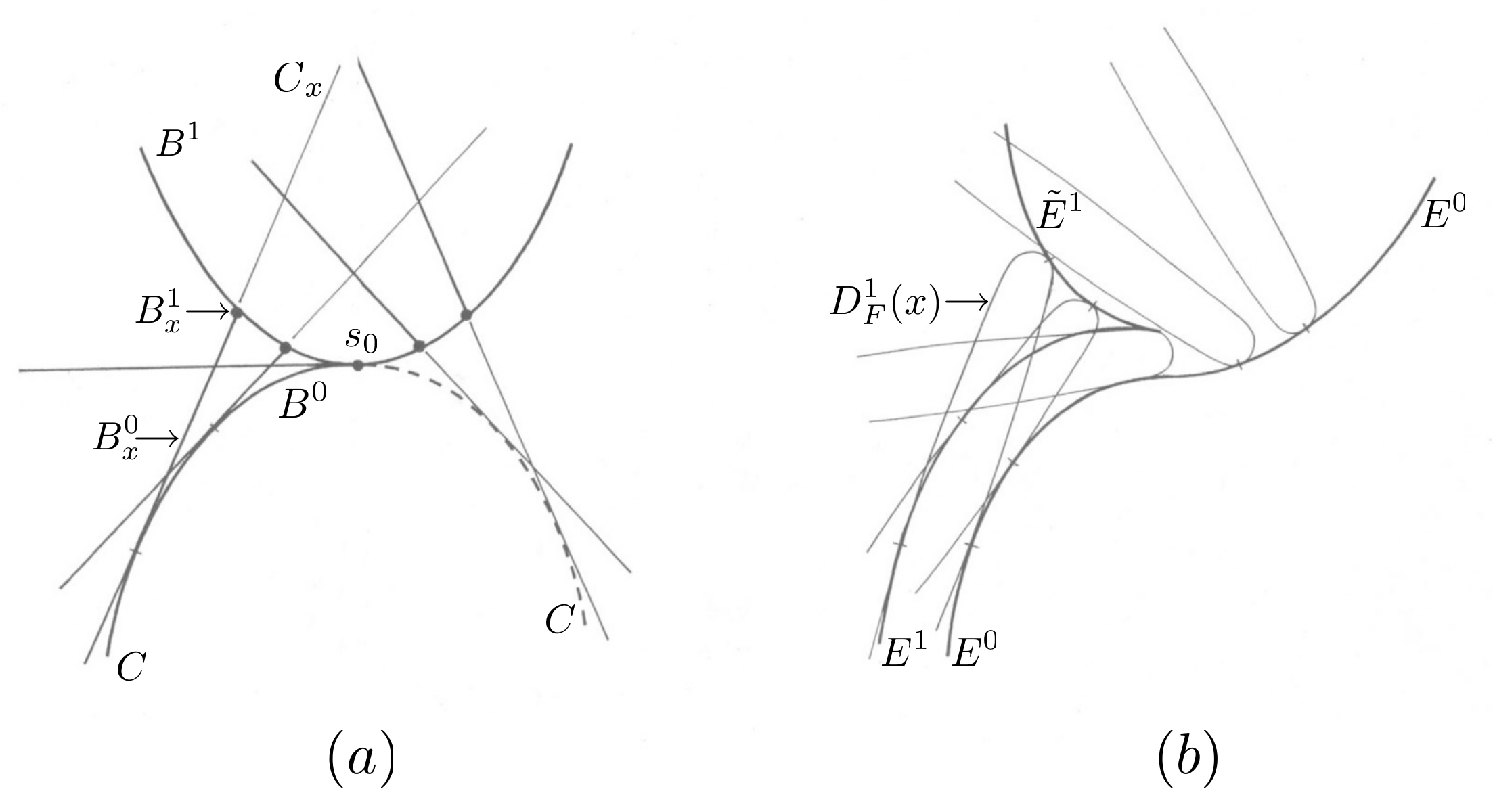} 
   \caption{Envelopes in $\bs^2$, sketched in the plane:
     (a) Local configuration of the envelope $B^0$ of arcs
     $B_x^0$ and the locus $B^1$ of their endpoints, 
     (b) corresponding configuration of envelope of smooth curves $D_F^1(x)$.} 
   \label{fig:ZZ}
\end{figure} 

We could leave the discussion at this point, but further useful insight is
obtained by taking account of the role of $b^{-1}(T)$, which we regard
(compare Section~\ref{s:d=k=1}) as a `first approximation' to $D_F^1(x)$.    
Let $C$ be the envelope of the great circles $\{C_x\}_{x\in S}$ and 
$B^0\subset C$ be the envelope of the semi-circular (open) arcs
$\bigl\{B_x^0\bigr\}_{x\in S}$. Let $B^1$ be the locus of
end-points $\bigl\{B_x^1\bigr\}$ for $x\in S$.  It is
straightforward to check that if $s_0\in C\cap B^1$ and is a fold
point of $C$ then $B^1$ is a smooth curve quadratically tangent to $C$
at $s_0$ where an arc of $B^0$ terminates: see Figure~\ref{fig:ZZ}(a). 
When we now replace $B_x^0$ by the figure-eight $D_F^1(x)$ then the
`approximate' envelope $B^0\cup B^1$ becomes a generic envelope
configuration $E$ consisting of a fold curve $E^0$ close to $B^0\,$ 
which continues as a fold curve close to one branch of $B^1$, 
together with another fold curve $E^1$ close to $B^0$ which meets a fold curve
$\til E^1$ close to the other branch of $B^1$ at a cusp point close to $s_0\,$:
see Figure~\ref{fig:ZZ}(b).  The numbers $0,2,4$ indicate relative
numbers $n_s$ for $s$ in each of the three complementary regions of
$E$ locally.  Finally, local stability implies that this configuration
occurs in $\rho\bs^2\cap\Delta_F$ for $\rho$ sufficiently small.
\smallskip

Features of the earlier $q=2$ case can be seen by intersecting the
above with a typical great circle $\bs^1$: we find interior points in
pairs (corresponding to the pair of branches of $E$ close to $C_0$)
and isolated end points (corresponding to the branches of $E$ close to
$C_1\,$).  The $q=3$ configuration shows how in the presence of a
further bifurcation parameter these features generically coalesce.  
\smallskip

Returning to the exceptional points $\{x_j\}$ on $S$ where $\rank
D_j=2$, let us make the generic assumption that for each $x_j$ neither
$1$-dimensional stratum $T_1$ or $T_1'$ of $T$ lies in the range of $D_j$.
It follows that for $x\in S$ close to $x_j$ the image $b_x(\bs^2)$ is
an oblate ellipsoid with `equatorial' plane close to the range of $D_j$.  
From this and Figure~\ref{fig:BB} it can be seen that $D_F^1(x)$ has 
the figure-eight configuration with one loop much larger than the other. 
As $x$ passes through $x_j$ the
smaller loop vanishes (the end point coalesces with the intersection point)
and is re-created at the antipodal end point; the transition curve $D_F^1(x_j)$
in $\bs^2$ has a $\frac32$-power cusp at both its antipodal end points.

Provided $x_j$ does not correspond to a cusp point of the envelope $E$
the local fold structure of the envelope itself (that is, the
bifurcation set in $\bs^2$) is unaffected as $x$ passes through $x_j\,$, 
although the location of corresponding solutions to the bifurcation
problem is controlled by this extra twist.  The bifurcation structure
at $x_j$ corresponds to one of the generic singularities of
$1$-parameter families of plane curves at the envelope: see~\cite{DU}. 
\medskip

For $q\ge4$ generically for all $x\in S$ the map $a({\cdot},x):\br^q\to\br^3$ is
a submersion at the origin: thus the local structure of $\Delta_F$ is
a product $\br^{q-3}\times\Delta_F^3$ where $\Delta_F^3$ is as in the
case $q=3$.  Nevertheless, the envelope of the family
$\bigl\{D_F^1(x)\bigr\}_{x\in S}$ in $\bs^{q-1}$ may have additional
singular structure, corresponding to higher-dimensional
generalisations of Figure~\ref{fig:ZZ} which we do not investigate here. 
\bigskip

We now turn to cases of greater degeneracy degree $m>2$.  
The essential geometry, while in principle more complicated,
has much in common with the case $m=2$ above.
\subsection{The case $m>2$.} \label{ss:m>2}  
Here
$$F(0,x,y)=h(y)=(0,y^m)$$
and this germ has $\K$-codimension $2m-1$ (see \cite{AVG}, \cite{GI}) 
with $\K$-miniversal deformation (using more systematic
notation than previously) of the form
\begin{equation}  \label{e:Hunf}
H(\til a,y)=(a_0+a_1y+\cdots+a_{m-1}y^{m-1},\,
                     \ba_0+\ba_1y+\cdots+\ba_{m-2}y^{m-2}+y^m)
\end{equation}
where 
$\til a=(a,\ba)=(a_0,\ldots,a_{m-1},\ba_0,\ldots,\ba_{m-2})\in\br^{2m-1}$.
The discriminant $\Delta_H$ is given by $R_m=0$ where $R_m$ is the
resultant of the two components of $H(\til a,y)$ and which has 
the following expression as a determinant
that for clarity we illustrate for $m=4$, 
the general pattern being analogous:
\begin{equation}  \label{e:r4}
R_4=
\left| \begin{array}{lllllll}
a_0 & a_1 & a_2 & a_3 & {} & {} & {} \\
{}  & a_0 & a_1 & a_2 & a_3 & {} & {} \\
{}  & {} & a_0 & a_1 & a_2 & a_3 & {} \\
{}  & {} & {} & a_0 & a_1 & a_2 & a_3 \\
\ba_0 & \ba_1 & \ba_2 & 0 & 1 & {} & {} \\
{} & \ba_0 & \ba_1 & \ba_2 & 0 & 1 & {} \\
{} & {} & \ba_0 & \ba_1 & \ba_2 & 0 & 1
\end{array} \right|.
\end{equation}
Here the blanks represent zero entries.
\begin{prop}
The resultant $R_m$ has the form
\begin{equation}  \label{e:res}
R_m=a_0^m + a_0P_m + (-1)^m\ba_0(a_1^m+a_1Q_m) + \ba_0^2(a_2^m+ \til R_m)
\end{equation}
where $P_m,Q_m$ are polynomials each of least degree $m$ and $\til R_m$ is 
a polynomial of degree at least $m+1$. 
\end{prop}
Here the polynomials are in the variables $a_i,\ba_j$ and \emph{least
degree} means the degree of the lowest-order term.\\[3mm]
\proof
The lowest-order term is $a_0^m$ from the main diagonal. Putting $a_0=0$
we then have $R_m=(-1)^{m}\ba_0 R_{m-1}'$ where $R_{m-1}'$ is a determinant
of order $2m-2$ illustrated in the case $m=4$ as follows:
\begin{equation}
R_3'=
\left|
\begin{array}{llllll}
a_1 & a_2 & a_3 & {} & {} & {} \\
{}  & a_1 & a_2 & a_3 & {} & {} \\
{}  & {} & a_1 & a_2 & a_3 & {} \\
{}  & {} & {} & a_1 & a_2 & a_3 \\
\ba_0 & \ba_1 & \ba_2 & 0 & 1 & {} \\
{} & \ba_0 & \ba_1 & \ba_2 & 0 & 1 \\
\end{array}
\right|.
\end{equation}
The lowest-order term in $R_{m-1}'$ is $a_1^m$, and then if $a_1=0$
the $\ba_0$ in the first column is multiplied by
a determinant whose lowest-order term is $(-1)^ma_2^m$.
\endproof
\medskip

Since the term of lowest order in $R_m$ is $a_0^m$ we see immediately:
\begin{cor}
The positive tangent cone $T_H$ lies in the hyperplane $a_0=0$.
\qed
\end{cor}
As in Section~\ref{ss:pull-back} we take polar coordinates 
$\eps=\rho s$ and write
$$(a,\ba)=(a(\rho s,x),\ba(\rho s,x))=\rho(b(s,x),\bb(s,x))+\rho^2(c(s,x),
                \bar c(s,x))+O(\rho^3)$$
where $(b,\bar b):\bs^{q-1}\times S\to\br^{2m-1}$ and likewise 
$(c,\bar c)$.  We then find from~(\ref{e:res}) that 
$R_m=\rho^{m}\bar{R}_m$ where
\begin{equation}  \label{e:rbar}
\bar R_m = \bar R_m(\rho,b,\bar b)
         =b_0^m + \rho\bigl((-1)^m\bb_0b_1^m + b_0\til P_m\bigr) 
                              + \rho^2\til Q_{m+2} + O(\rho^3)
\end{equation}
in which $\til P_m,\til Q_{m+2}$ are homogeneous polynomials of degree $m,m+2$
respectively in the variables $(b_i,\bb_j)$.
Following the same procedure as in the case $m=2$ we first study the
geometry of the set
\[
\Delta_F^1:=\{(\rho,s,x)\in\br\times\bs^{q-1}\times S:
                      \bar R_m(\rho,b(s,x),\bar b(s,x))=0\}
\]
since $\bar R_m(\rho,b,\bar b)=0$ for $\rho\ne0$ precisely when
$(a,\bar a)\in\Delta_H\,$,
and then apply the pinching map $\til\tau$ to obtain 
$\Delta_F=\til\tau(\Delta_F^1)\subset\br^q\times S$.
\medskip

We shall require some transversality hypotheses analogous to (H1).
Rather than formulate a general statement (which we do not need and in
any case cannot carry out without specifying a stratification of $T_H$) we
state the appropriate hypotheses as we come to them; they mirror those of
Section~\ref{sss:geom}.
\medskip

Suppose $(s_0,x_0)$ satisfies $b_0(s_0,x_0)=0$, and take a local
coordinate chart $(s,x)=(s_0+v,x_0+w)$ on $\bs^{q-1}\times S$ as before.  
As in the $m=2$ analysis we look at three cases.
\bigskip

(\emph{i})\quad $\bb_0b_1\ne0$ at $(s_0,x_0).$
\medskip

Here by the IFT the locus $\bar R_m=0$ locally has the form 
\begin{equation}  \label{e:regcase}
\rho=k_0 b_0^m + l_{m+1}\,,\quad k_0\ne0,
\end{equation}
where $b_0=b_0(s_0+v,x_0+w)$ and where $l_{m+1}$ denotes terms of
degree at least $m+1$ in $v,w$.
We now make the generic hypothesis 
\medskip

\noindent $\bigl($H$1_m(i)\bigr)$ \quad $(s_0,x_0)$ is a regular point 
of the function $b_0$.
\medskip

Here $\Delta_F^1$ is locally a smooth hypersurface in
$\br\times\bs^{q-1}\times S$ having contact of order $m$ with 
$0\times\bs^{q-1}\times S$ along the codimension-$1$ submanifold
$$B_0:=\{(s,x)\in \bs^{q-1}\times S: b_0(s,x)=0\}.$$
Applying the pinching map $\til\tau$, remembering that
$\tau(\rho,-s)=\tau(-\rho,s)$ and using Proposition~\ref{p:pinch}
we obtain $\Delta_F$ locally as a
$C^1$ hypersurface in $\br^q\times S$ containing $0\times S$ and having contact
of order $(m+1)/m$ along $0\times S$ with the \lq direction locus' in
$\br^q\times S$ determined by $B_0\,$.  
Thus for $q=2\,$, when $m=2$ the discriminant $\Delta_F$ has a
$\frac32$-power cusped ridge along $0\times S$ as already seen in
Section~\ref{sss:q=2}, while for $m=3$ it is a $C^1$ submanifold of
$\br^2\times S$ passing through $0\times S$. 
\bigskip

(\emph{ii})\quad $\bb_0=0, b_1\ne0$ at $(s_0,x_0)$.
\medskip

Here we make the next generic hypothesis:
\bigskip

\noindent$\bigl($H$1_m(ii) \bigr)$\quad $(s_0,x_0)$ is a regular zero of the map 
$(b_0,\bb_0):\bs^{q-1}\times S \to\br^2$.
\bigskip

First take $q=2$.
The transversality hypothesis (H$1_m(ii)$) means that we can regard the pair
of functions $(b_0,\bb_0)$ as local coordinates on $\bs^1\times S$ close
to $(s_0,x_0)$.   We can then write $\bar R_m$ in the form
\begin{equation}  \label{e:pform}
p^m+\rho\bigl(\bar p + c\rho + O(2)\bigr)
\end{equation}
where $(p,\bar p):= \bigl(b_0,(-b_1)^m \bb_0 + b_0\tilde P_m\bigr)$ and
$c$ is a constant.  The quadratic terms are 
nondegenerate regardless of~$c\,$, so by the Splitting Lemma 
(see~\cite{CH},\cite[\S4.5]{Dem},\cite{POS}) we can suppose the coordinates
are such that $\Delta_F^1$ has the form
\begin{equation} \label{e:endlocus}
p^m+\rho(\bar p + c\rho)=0.      
\end{equation} 
The description of this locus depends on whether $m$ is even or odd.
\smallskip

The case when $m$ is \emph{even} is relatively easy.  If
$m=2l$ then we write $t=p^l$ and the left hand side 
of~(\ref{e:endlocus}) becomes a nondegenerate quadratic form in
$(\rho,t,\bar p)$: the zero set is a cone with vertex at the origin,
tangent
to the $(t,\bar p)$-plane and with the $\bar p$-axis as a generator.  On
reverting to $(\rho,p,\bar p)$-coordinates the cone retains these features,
although the contact with the plane $\rho=0$ is now of
order $m$ rather than quadratic.  
Applying $\til\tau$ then gives $\Delta_F$ with a \lq ship's prow' geometry
at $(0,x_0)$ in $\br^q\times S$ as in Section~\ref{sss:q=2} (see
Figure \ref{fig:DD}) but with cusp contact of order $(m+1)/m$.
\smallskip

\begin{figure}[htbp] 
   \centering
   \includegraphics[scale=.5]{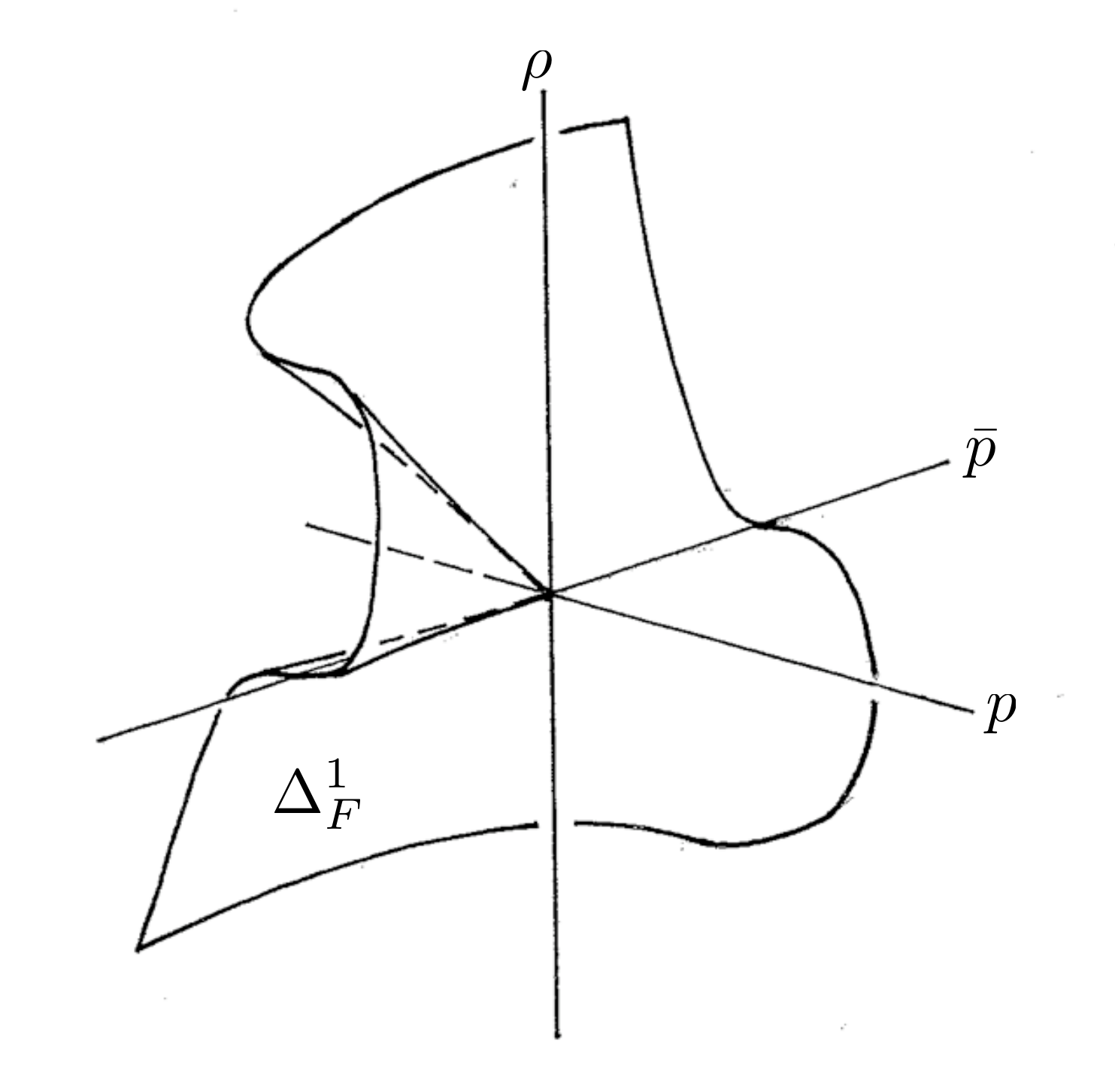} 
   \caption{The solution locus $\Delta_F^1$ given by  
   $\overline{R}_m=0$ for $m=3$, $q=2$  at a point where 
   $b_0=\overline{b}_0=0$. Here the picture is constructed by taking   
   $\overline{R}_3=p^3+\bar{p}\,\rho+\rho^2$ where $p,\bar{p}$ are the 
   coordinates $p=b_0,~~\bar{p}=(-1)^mb_1^m\,\bar{b}_0+b_0\tilde{P}_m$
   (see text).}
   \label{FigureXZ}
\end{figure}

The case when $m$ is \emph{odd} is a little harder to describe pictorially.
Some insight is gained by regarding~(\ref{e:endlocus}) as the equation
of a curve in the $(\rho,p)$-plane for each fixed $\bar p$.
As $\bar p$ increases through zero, the order-$m$ tangency of the locus with
the plane $\rho=0$ along the $\bar p$-axis reverses its orientation as
a `bulge' passes through the origin from positive to negative $\rho$
where as before we interpret $(-\rho)s$ as $\rho(-s)$ when $\rho<0$.
A 3-dimensional representation is attempted in Figure \ref{FigureXZ}.
%
\smallskip

For the general case when $q>2$ the local
structure of $\Delta_F$ is the cartesian product of the
above with $\br^{q-2}$.
\bigskip

(\emph{iii})\quad $\bb_0\ne0, b_1=0$ at $(s_0,x_0)$.
\medskip

Substituting $a_0=a_1=0$ into the expression~(\ref{e:res}) for $R_m$
we obtain $R_m=\ba_0^2D_m$ where $D_m$ is a polynomial whose term of
lowest order is $a_2^m$.   Therefore with $q=2$ we make the further
generic hypothesis
\bigskip

\noindent $\bigl($H$1_m(iii)\bigr)$\quad $b_2(s_0,x_0)\ne0$.
\bigskip

It then follows that $(\eps,x)\notin\Delta_F$ for sufficiently small 
$|x-x_0|$ and $|\eps|$.
Thus there is generically no analogue of intersection points in 
$\Delta_F$ when $m>2$.  
\smallskip

When $q>2$ there exist generically points at
which $b_0=b_1=b_2=0$, which we do not investigate further here.
\subsubsection{The bifurcation set for $m>2$.}  \label{ss} 
Having found the local form of the discriminant
$\Delta_F\subset\br^q\times S$ under
generic assumptions, we are now in a position to describe the
generic structure of the bifurcation set, that is the apparent outline 
of $\Delta_F$  viewed along $S$. Recall that $\Delta_F$ is the locus 
$\overline R_m=0$ where
$\overline R_m$ is given by~(\ref{e:rbar}) with $(a,\ba)=\rho(b,\bb)$. 
Again we work first with $q=2$.
\medskip

Taking local coordinates $(s,x)=(s_0+v,x_0+w)$ on $\bs^1\times S$ we
seek common solutions to the equations
\[
\overline R_m=0\,,\qquad \dbd{}w\overline R_m=0.
\]
As before, we consider two cases.
\bigskip

(\emph{i})\quad $\bb_0b_1\ne0$ at $(s_0,x_0)$.
\medskip

Here the argument parallels that of Section~\ref{sss:bifgeom}. We make
the generic hypothesis
\bigskip

\noindent $\bigl($H$2_m(i)\bigr)$\quad The projection 
$\mu|B_0:B_0\to\bs^1$ has only {\em fold} singularities. 
\bigskip

If $(s_0,x_0)$ is a regular point of $\mu|B_0$ then no
local bifurcations occur for any solution branches in directions
sufficiently close to $s_0\,$.  If $(s_0,x_0)$ is a fold point of
$\mu|B_0$ then from~(\ref{e:regcase}) and following Section~\ref{sss:bifgeom} 
we find that the apparent outline of $\Delta_F^1$ in the direction of
$S$ is locally a smooth curve $\Gamma_0$ through $(0,s_0)$ and having
contact of order $m$ with the $(v,w)$-plane, so that in
particular when $m$ is odd it has only one branch for $\rho>0$ while for
$m$ even it has two.  Applying the pinching map $\til\tau$ to
$\Delta_F^1$ and projection $\mu$ gives a curve $C_0$ through the origin in
parameter space $\br^2$ having contact of order $(m+1)/m$ with its tangent direction
$s_0\,$.  The number of solutions $(x,y)$ changes by two as
each branch of $C_0$ is crossed.
\bigskip

(\emph{ii})\quad  $\bb_0=0$,$\,b_1\ne0$ at $(s_0,x_0)$.
\medskip

Here we make the further generic hypothesis 
\bigskip

\noindent $\bigl($H$2_m(ii)\bigr)$\quad 
$\dbd{b_0}x\bigl((-b_1)^m\dbd{\bar b_0}x + \til P_m \dbd{b_0}x\bigr) \ne 0$
at $(s_0,x_0)\,$. 
\bigskip

\begin{figure}[htbp] 
   \centering
   \includegraphics[scale=.5]{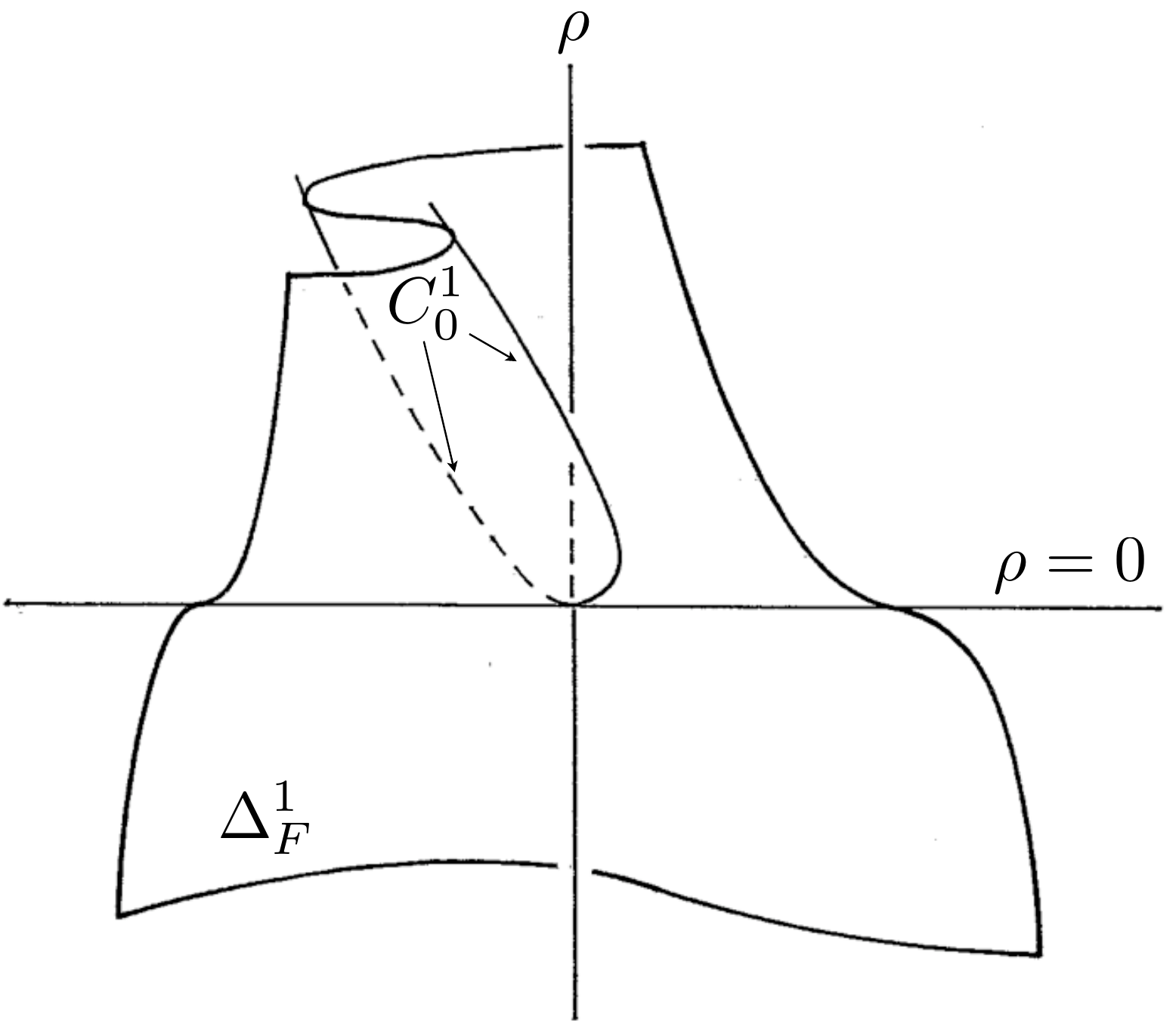} 
   \caption{The discriminant $\Delta_F^1$ for $m=3$ viewed in a
     typical direction not parallel to $p,\bar{p}$ axes: apparent outline
     ${C}^1_0$ (bifurcation curve) in $\br^2$ becomes a cusped 
     curve at the origin after applying the pinching map $\tau$.} 
     \label{fig:YZ}
\end{figure}

In terms of coordinates $(p,\bar p)$ used above, this simply states
that $\xi:=\dbd px$ and $\eta:=\dbd{\bar p}x$ are both nonzero 
at $(p,\bar p)=(0,0)$.  From~(\ref{e:pform}) and again
disregarding the higher order terms as in~(\ref{e:endlocus})
we have 
$$\dbd{\overline R_m}x = mp^{m-1}\xi + \rho\eta$$
so that $\dbd{\overline R_m}x = 0$ is solved as
$\rho=\til\rho(p)=\lambda p^{m-1}$ with constant
$\lambda=-m\xi/\eta\ne0$ at $(0,0)$.  
Substituting into~(\ref{e:endlocus}) and deleting the factor $p^{m-1}$ we obtain
\begin{equation}  \label{e:barb}
p+\lambda\bar p +O(p^{m-1})=0
\end{equation}
which represents a smooth curve $B_0'$ through the origin in the
$(p,\bar p)$-plane, transverse to both axes and not tangent to
$(\xi,\eta)$ since $m\ne1$. Projecting the graph of
$\til\rho|B_0'$ to a plane orthogonal to $(\xi,\eta,0)$ then 
gives the apparent outline of $\Delta_F^1$ viewed along $S\,$: 
this curve $C_0^1$ has one branch for
$\rho>0$ and one for $\rho<0$ when $m$ is even and two branches for
$\rho>0$ when $m$ is odd. After applying the pinching map $\tau$ the curve
$C_0^1$ becomes a bifurcation curve $C_0$ having contact of order
$(m+1)/m$ with its tangent line at the origin.
For $m$ even the number of solutions
changes locally from $0$ to $2$ as the curve is crossed; for $m$ odd it
changes from $1$ to $3$ and back again in a hysteresis loop as the two 
branches are successively crossed.  See Figure~\ref{fig:YZ}.
%
\medskip

We summarise the bifurcation results in the following theorem.  Recall
that we are in the same dimensional setting as Theorem~\ref{t:bif1}
but now with degeneracy degree $m>2$.

\begin{figure}[htbp] 
   \centering
   \includegraphics[scale=0.7]{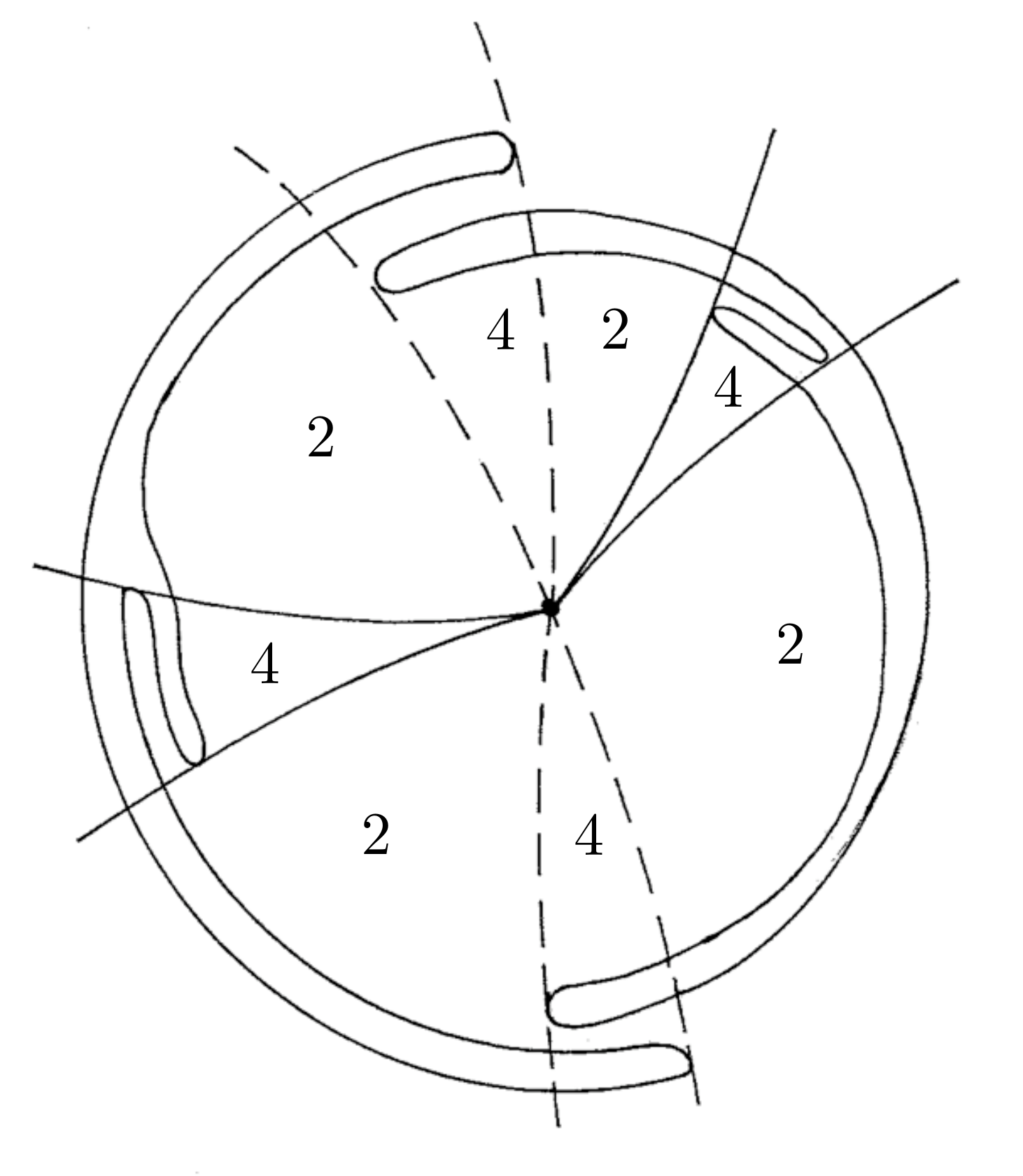} 
   \caption{Schematic bifurcation diagram  for $d=k=1, q=2$ and $m$
     odd (see text), showing hysteresis arcs (unbroken lines) and fold
     arcs (dashed lines) through the origin in parameter space
     $\br^2$. Numbers indicate the number of solutions in
     each region. Note that for each opposite pair of fold arcs the
     two folds approach from opposite sides in the plane.}
   \label{fig:figbifB}
\end{figure}

\begin{theo}  \label{t:bif2}
Under the generic hypotheses {\em $\bigl($H$1_m(i)-(iii)\bigr)$} and 
{\em $\bigl($H$2_m(i)-(ii)\bigr)$}
the bifurcation set for~(\ref{e:orig}) consists of a finite set of 
$C^1$ arcs emanating from the origin
in $\br^2$, each making contact of order $\frac{m+1}m$ 
with its tangent ray at the origin.  The arcs are of two types:
\emph{fold} and \emph{end} arcs, corresponding to the
following bifurcation behaviour as $\eps$ rotates about the origin 
in $\br^2$ in (locally) the appropriate direction:
\medskip

\emph{Fold arcs}.\quad 
These occur in pairs with a common tangent at
the origin, emanating in the same or opposite directions from the
origin according as $m$ is even or odd respectively. In the \emph{even} case,
as $\eps$ crosses the first arc in the appropriate direction
a pair of solutions is created away from $S$ by a saddle-node (fold)
bifurcation, and as $\eps$ crosses the second arc 
a second saddle-node pair of solutions is
created close to the $\kappa$-images of the first pair.  
In the \emph{odd} case a saddle-node bifurcation occurs across each arc.
\medskip

\emph{End arcs ($m$ even)}.\quad 
These occur in pairs with a common tangent at
the origin, emanating in opposite directions from the
origin.  As $\eps$ crosses an end arc in the appropriate
direction a pair of solutions is created away from $S$ in a 
saddle-node bifurcation; the two directions are in opposing senses
of rotation about the origin.
\smallskip

\emph{Hysteresis arcs ($m$ odd)}.\quad
These occur in pairs with a common tangent at
the origin, emanating in the same direction from the
origin.   As $\eps$ crosses the first arc in the appropriate
direction a saddle-node bifurcation creates a pair of
solutions, one of which coalesces with an existing solution 
(via a hysteresis loop) in a saddle-node bifurcation
across the second arc.   
\end{theo}
This elaborates on and gives a geometric setting for Theorem~2.1B
of~\cite{HT2}. See Figure~\ref{fig:figbifB}, and compare with
Figures~3 and~4 of~\cite{HT2}.
%
Note that the interpretation regarding solutions $(x,y)$ passing
through $y=0$ does not hold in the general case $m>2$ since 
$\bar a_0=0$ does not imply $a_0=0$ for points on $\Delta_H\,$. 
\medskip

\rem These results on the structure of $\Delta_F$ and the form of the 
expression~(\ref{e:Hunf}) show that for $\eps\in\br$ and $d=k=1$ the
generic solution branches are of the form $(\eps,x,y)\sim(t^m,t^{m+1},t)$.  
This applies also in the case $m=1$ (normal nondegeneracy) where $y$
represents a coordinate transverse to $S$ and the nonzero solution
value $y$ varies linearly with $\eps$. 
\subsection{The variational case}   \label{ss:varcase}
As noted in Section~\ref{ss:gradient}, 
when the original vector field $F_\eps$ is the gradient of a smooth real-valued
function $f_\eps$ then 
the reduced vector field $\til F_\eps$ may be regarded
as the gradient of a real-valued function $\til f_\eps\,$.  Therefore
in this setting the bifurcations of equilibria are controlled not as
in Section~\ref{s:versal} by $\K$-versal deformation of the map germ
\[
\til F_0(x_0,{\cdot}):\br^k\to\br^d\times\br^k
\]
at $0\in\br^k$ but rather by an 
${\mathcal R}$-versal deformation of the function germ  
\[
\til f_0(x_0,{\cdot}):\br^k\to\br
\]
at $0\in\br^k$ where $\til F_0=\grad \til f_0 $.  
In the case of uniform normal degeneracy of $S$ the
germ $\til f_0$ is independent of $x_0\in S$. 
\medskip

 To study the bifurcation
geometry we follow the formal procedure of Section~\ref{ss:defman} but
now in addition to the simplification in dealing with functions rather
than vector fields there is an added twist: the gradient of a
deformation of $\til f_0(x_0,{\cdot})$ will involve derivatives with
respect to the `parameter' $x$ as well as with respect to the variable $y$.
\medskip

Assume we are in the uniform case, and write $\til f_0(x,y)=\til f_0(y)$.
If $\til f_0$ has an $r$-parameter $\R$-versal deformation of the form 
\[
a_1v_1(y)+\cdots+a_rv_r(y)+\til f_0(y)
\]
then as in Corollary~\ref{c:global} we express $\til f_\eps$ as
\[
\til f_\eps=a_1(\eps,x)v_1(\til y) +\cdots+ a_r(\eps,x)v_r(\til y) +
\til f_0(\til y) 
\]
and now (after dropping the tilde)
we find 
$F_\eps=(\partial_x f_\eps,\partial_y f_\eps)$ 
where
\begin{equation}
\begin{array}{ll}
\partial_x f_\eps = a_1'(\eps,x)v_1(y)+\cdots+a_r'(\eps,x)v_r(y)\\[3mm]
\partial_y f_\eps = a_1(\eps,x)v_1'(y)+\cdots+a_r(\eps,x)v_r'(y)
                                       +F_0(y)
\end{array}
\end{equation}
and we write $a_i',v_j'$ for
$\grad_x a_i$, $\grad_y v_j$ respectively.
The condition for $(\eps,x)$ to belong to the
discriminant $\Delta_F$ is an algebraic condition on the $1$-jet 
(with respect to $x$) of the map $a=(a_1,\ldots,a_r):\br^q\times
S\to\br^r$.  
We illustrate this with a simple but non-trivial example that exhibits
both the similarities with and the differences from the general
geometric formalism of previous sections.
\bigskip

\textbf{Example.}\quad  Let $k=1$ and consider the bifurcation problem
$F_\eps=0$ where $F_\eps=\grad f_\eps$ with
\[
f_\eps(x,y)=\eps g(\eps,x,y)+y^m.
\]
An $\mathcal R$-versal deformation of the function germ $y^m$ is 
\[
a_0+a_1y+\cdots+a_{m-2}y^{m-2}+y^m
\]
(see~\cite{AVG},\cite{GI}) and so for $m=3$ we can write
\[
f_\eps(x,y)=a_0(\eps,x)+a_1(\eps,x)y+y^3
\]
and we have
\[\begin{array}{ll}
\partial_xf_\eps(x,y)=a_0'(\eps,x)+a_1'(\eps,x)y  \\[3mm]
\partial_yf_\eps(x,y)=a_1(\eps,x)+3y^2.
\end{array}\]
The condition for $(\eps,x)$ to belong to the discriminant $\Delta_F$
in the case $d=1$ (compare~(\ref{e:delH})) is then
\[
3(a_0'(\eps,x))^2+(a_1'(\eps,x))^2a_1(\eps,x)=0.
\]
Following the geometric analysis of Section~\ref{s:d=k=1} we write
$\eps=\rho s$ in polar coordinates and 
\begin{equation} \label{e:rhob}
a(\rho s,x)=\rho b(s,x)+O(\rho^2)
\end{equation}
and find that the conditions for $b(\eps,x)$ to lie in the positive
tangent cone $T$ are
\begin{equation}
b_0'(x)=0, \quad b_1(s,x)\le 0\,.
\end{equation}
where $b_0'=\grad_xb_0$.
For simplicity of description we take as usual the case $q=2$ so $(s,x)\in
\bs^1\times S$. Generically the locus 
$$
B_0'=\{(s,x)\in \bs^1\times S:b_0'(s,x)=0\}
$$ 
is a collection of smooth curves (codimension-$1$ manifold) 
such that the singularities of the projection
$\mu|B_0':B_0'\to \bs^1$ are isolated folds.  There are at least two
of these on each compact connected component of $S\,$.  
(For $q\le4$ the singularities of $\mu$ will generically be folds, cusps ($q=3,4$) or
swallowtails ($q=4$), while for $q>4$ singularities of corank $>1$ can
appear.  See~\cite[Ch.VII,\S6]{GG} for example.)
Likewise the locus
$$
B_1=\{(s,x)\in \bs^1\times S:b_1(s,x)=0\}
$$
is generically a smooth curve (codimension-$1$ manifold) which
intersects $B_0'$ transversely at points which are \emph{not} singular
points of the projection $\mu|B_0'$.   
The analysis proceeds as in Sections~\ref{sss:q=2} and \ref{sss:bifgeom}
noting that now {\em end-points} are given by $B_0'\cap B_1$
while {\em intersection points} are the critical points of
$b_1|B_0\,$. Thus the generic bifurcation geometry turns out to be the same
as in Theorem~\ref{t:bif1} for the (non-variational) case $m=2$
although the interpretation of key features is here tied to the
variational structure.
\bigskip

For the general case $f_0(y)=y^m$ we define the $1$-jet map
\[
j^1a:\br^q\times S \to  J^1(S,\br^{m-1}):(\eps,x)\mapsto 
              (x;a_0,a_1,a_2,\ldots,a_{m-2};a_0',a_1',\ldots,a_{m-2}')
\]
where $a_i=a_i(\eps,x)$ and $a_i'=\grad_xa_i$ for $i=0,1,\ldots,m-2$.
Taking local coordinates on $S$ we consider the map
\[
h:\br\to\br^{d+1}:y\mapsto (0,\ldots,0,my^{m-1})
\]
which has $\K$-codimension $e:=(d+1)(m-1)-1$ (see~\cite{GI}) and a $\K$-miniversal
deformation $H$ with discriminant $\Delta_H\subset\br^e$. Locally
$J^1(S,\br^{m-1})\cong \br^d\times\br^{m-1}\times\br^{d(m-1)}
=\br^d\times\br^{e+1}$ and we
let $\Delta_f=(\pi\circ j_1a)^{-1}\Delta_H$ where $\pi$ is projection
of $J^1(S,\br^{m-1})$ onto the last $e$ coordinates (that is, omitting
$x$ and $a_0$). We then consider the related map 
\[
\til b:\bs^{q-1}\times S\to\br^e:(s,x)\mapsto 
              (b_1,b_2,\ldots,b_{m-2},b_0',b_1',\ldots,b_{m-2}')
\]
with $b:\bs^{q-1}\times S\to\br^{m-1}$ defined as in~(\ref{e:rhob}).
Under suitable non-degeneracy hypotheses on $\til b$ which, using the
Thom jet transversality theorem
(see~\cite{AVG},\cite{Dem},\cite{GI},\cite{GG}), we expect
to be generically satisfied by the map $f_\eps\,$, local geometry of
the discriminant $\Delta_f$ and its projection to the parameter space $\br^q$
(giving the bifurcation structure) can in principle be determined --- at least for
small values of $q$ and $m$.  We leave details to the enthusiastic reader.  
\subsection{Possible further reduction}
In the general case a further reduction is technically  
possible, using the IFT one more time.  
For simplicity of illustration take $q=1$ so that $\eps\in\br$.
\medskip

Given $x\in S$, write $\til g(x)=\tilF(0,x,0)$ where $\tilF:\br\times
U\to\br^{d+k}$ is the deformation map as in~(\ref{e:calF}) with
$(x,y)\in U$.  We make the following assumption, which holds
generically as $\dim S=d<d+k\,$:
\bigskip

\noindent (A4)\quad
 The image of $\til g:S\to\br^{d+k}$ avoids the origin in $\br^{d+k}$.
\bigskip

\noindent 
This means that given $x_0\in S$ there is at least one component $\til
g_i$ of $\til g$ that does not vanish at $(x_0,0)$ when $\eps=0$.  By the IFT
the $i^{\mbox{th}}$ component of $\til F_\eps=0$ can then be solved
close to $u=(x_0,0)$ in the form $\eps=\eps_i(x,y)\,$.  Substituting
this into the other components of $\til F_\eps$ yields $d+k-1$
equations to be solved for $(x,y)$ close to
$(0,0)\in\br^d\times\br^k$.  
The solution set can be viewed as
the zero set for a $d$-parameter family of smooth map germs
$\br^k\to\br^{d+k-1}$, where techniques of versal deformation still
apply, now in one lower dimension in parameter space and target.
\medskip

For example, if $d=1$ then we have a $1$-parameter deformation of
germs $\br^k\to\br^k$ about which much is known for small $k$:
see~\cite{RI} for the case $k=2$.  If $d=2$ then we are dealing with
germs $\br^k\to\br^{k+1}$: see~\cite{MD} for the case $k=2$ here.
\medskip

In this formal reduction of the problem, however, the controlling role
of $\eps$ has been sacrificed, and it is not clear that there is
practical benefit to understanding the bifurcation structure.  
We do not pursue this approach further.
\subsection{Nontrivial bundles}  \label{ss:nontrivial}
In the case of nontrivial bundles $NS,TS,K$ etc, where the same
coordinates $(x,y)$ cannot be used on the whole of $U\subset K$, the
deformation theory must be re-expressed in terms of jet bundles.
Given smooth vector bundles $D,E$ over $S$ and a family of bundle maps
$F_\eps:D\to E$ with $\eps\in\br^q$, the generic structure of $F_\eps$
is expressed in terms of transversality of the $r$-jet extension
$$j^r F : \br^q\times D\to J^r(D,E): (\eps,u)\mapsto j^rF_\eps(u)$$ 
to suitable stratifications of the $r$-jet bundle $J^r(D,E)$ and for 
sufficiently large $r$ (cf.~\cite{AGLV}).  In our context
this theory has to be applied to the reduced vector field
$\til F_\eps :K\to TS\oplus P$.  
The geometry of the zero set is captured in local
coordinates using $\K$-versal deformation theory as we have done, but 
global aspects are dictated by the geometry of the bundles. 
It would be interesting to explore this further even in
the simplest case $d=k=1$ with $K$ both the trivial and nontrivial 
$\br$-bundle over a circle.  
\section{Local branching analysis in the case $q=k=1$.}  \label{s:q=k=1}
It is one thing to understand generic structure of the bifurcation
set, but another to recognise the specific branching geometry in
particular cases.  In this section we focus on the latter.
We return to the key case 
$q=k=1\,$, taking advantage of the fact that a singular function of one 
variable $y\in\br$ is dominated by the first nonvanishing term of its 
Taylor series in order to classify the local branching
behaviour of $\til F_\eps$ in a systematic algebraic way (bottom up)
that complements the generic geometrical approach (top down) that we have
described so far.    
\medskip

Recall from~(\ref{e:calF}) that we have in local coordinates
\begin{equation}  
\label{e:icomp}
\til F(\eps,u)=\til F_0(u)+\eps\tilF(\eps,u)
\end{equation}
where $u=(x,y)\in S\times\br$ and $\eps\in\br\,$.  Each component of
$\til F_0(u)$ has a Taylor series in $y$ with coefficients depending smoothly
on $x\in S$.  We do \emph{not} immediately assume here that
the normal degeneracy is uniform,  but we suppose some control on
degeneracy by making the following assumption:
\bigskip

\noindent (A5)\quad
Each component of $\til F_0(u)$ has nonzero Taylor series
in $y\,$ at every point $x\in S$.
\bigskip

In other words, the normal degeneracy in the problem is of finite order
everywhere on $S\,$.  This is not essential in all that follows, for it
will be clear how to amend the discussion to allow infinite degeneracy
in one component.   We leave aside, however, the complications of infinite
degeneracy in more than one component.
\subsection{Necessary conditions for branching} \label{ss:nec}
From (\ref{e:tilf}) and (\ref{e:dtilf}) 
the $i^{\mbox{th}}$ component of $\til F(\eps,u)$ can be written in the form
\begin{equation}  \label{e:Fi}
\til F_i(\eps,u)=\eps g_i(\eps,u)+y^{m_i}r_i(u) 
\end{equation}
for some integer $m_i\ge2$ and a smooth function $r_i:U\to\br\,$;
we suppose that $r_i({\cdot},0)$ is not identically zero, otherwise we
would increase the exponent $m_i\,$.  Let $x_0\in S$ and number the
coordinates so that at $x_0$
\[
m_1\le m_2\le\cdots\le m_{d+1}.
\]
\begin{prop}
\label{p:necessary}
Necessary conditions for $x_0\in S$ to be a branch point
are that for $1\le i < j \le d+1\,$
\begin{equation} 
\label{e:nconds}
\begin{array}{ll} 
r_ig_j=0\qquad \mbox{when}\quad m_i<m_j \\[3mm]
d_{ij}=0\qquad \mbox{when}\quad  m_i=m_j
\end{array}
\end{equation}
for all $1\le i < j \le d+1\,$, where $g_i=g_i(0,x_0,0)\,$, $r_j=r_j(x_0,0)$
and $d_{ij}:=r_ig_j-r_jg_i\,$.
\end{prop}
\proof
From (\ref{e:Fi}) we see that 
\begin{equation}  
\label{e:elim}
\til F_i(\eps,u)g_j(\eps,u)-\til F_j(\eps,u)g_i(\eps,u)
   = y^{m_i}r_i(u)g_j(\eps,u) - y^{m_j}r_j(u)g_i(\eps,u).
\end{equation}
Dividing through by $y^{m_i}$ and letting $(\eps,x,y)\to (0,x_0,0)$ 
immediately gives the results.
\endproof
\medskip

The conditions (\ref{e:nconds}) are not all independent.  
Given that the assumption (A5) holds, choose $x_0\in S$ and let 
$$n=n(x_0)=\max\{i:g_i\ne0\}\,.$$
Then the conditions (\ref{e:nconds}) reduce to
\begin{equation}
\label{e:nconds2}
\begin{array}{ll}
(a)\quad r_i=0\qquad \mbox{when}\quad m_i>m_n \\[3mm]
(b)\quad d_{ij}=0\qquad \mbox{when}\quad m_i=m_j=m_n\,.
\end{array}
\end{equation}
Let $J_n=\{j:m_j=m_n\}$: then $J_n$ is an integer interval $[l,p]$ say,
with $l\le n\le p$.  In~(\ref{e:nconds2}) the number of conditions $(a)$
is $l-1$ and the number of independent conditions $(b)$ is $p-l$.
Thus~(\ref{e:nconds2}) represents $p-1$ conditions, which together
with 
\[
g_{n+1}=\cdots=g_{d+1}=0
\]
give a total of $d-n+p$ conditions on $x_0\in S$.  Hence we have the
following result.
\begin{cor}  \label{c:isol}
Generically no points $x_0\in S$ satisfy the
conditions~(\ref{e:nconds2}) when $n<p$, and there is a set of
isolated points $x_0\in S$ which satisfy~(\ref{e:nconds2}) when $n=p$.
\end{cor}
In the uniform case the results can be expressed more simply.
Each $r_i$ is nonzero everywhere on $S$ and the exponents $m_i$ do not
depend on $x\in S\,$, so if branching is to occur at any $x_0\in S$ 
we must have $n=p$ (since if $n<i\le p$ then
$d_{ni}=0$ implies $r_i=0$) and also $l=1$ (else $r_j=0$ for
$j<l$).  Hence we obtain in this context the neater result:
\begin{cor}
In the uniform case a necessary condition for $x_0\in S$ to be a
branch point is that $g_{[n]}(x_0)=0$ where
\[
g_{[n]}:=(d_2,\ldots,d_n,g_{n+1},\ldots,g_{d+1}):S\to\br^d
\]
in which $d_j$ denotes $d_{1j}$ for $j=n,\ldots,d$ and 
$n=\max\{i:m_i=m_1\}$.
\end{cor}
Of course, if $m_1<m_2$ (so $m_1<m_j$ for all $j>1$) then $g_{[n]}$
contains no $d_j$ terms and the result is even more straightforward :
\begin{cor}
If $m_1 < m_2$ then branch points $x_0\in S$ must be
zeros of the map 
$$\bar g:=(g_2,\ldots,g_{d+1}):S\to\br^d.$$
\end{cor}
Generically $0\in\br^d$ will be a regular value of $\bar g$ and
the zeros of $\bar g$ will be a discrete subset of $S$. In any given
context there may be analytical or topological tools (for example,
index theory) available to locate or count them.
\bigskip

The previous results gave necessary conditions for solution
branches to emanate from $(0,x_0,0)\in\br\times U
\subset\br\times S\times\br$.  We now turn to sufficient conditions.
\subsection{Sufficient conditions for branching}
Let $x_0\in S$ satisfy the necessary conditions~(\ref{e:nconds2}) for
branching.  With the notation of Section~\ref{ss:nec}
since $g_n\ne0$ at $(0,x_0,0)$ we can use the IFT to solve the equation 
$\til F_n(\eps,u)=0$ from~(\ref{e:Fi}) in the form
\begin{equation} 
\label{e:meqn}
\eps=\eps_n(u)=y^{m_n}\del_n(u)
\end{equation}
where $\del_n$ is a smooth function on a \nhd of $(x_0,0)$ in
$U\subset K$.
Substituting this into each component $\til F_i$ with $i\ne n$ and
dividing through by appropriate powers of $y$ we see that finding a
solution branch amounts to finding solutions $u=(x,y)$ close to $(x_0,0)$ to
the set of equations:
\[\begin{array}{lll}
\del_n(u) g_i(\eps_n(u),u) + y^{m_i-m_n}r_i(u) =0,\quad i>n \\[3mm]
                d_{nj}(\eps_n(u),u) =0, \quad l\le j<n \\[3mm]
y^{m_n-m_j}\del_n(u)g_j(\eps_n(u),u) + r_j(u)  =0,\quad j<l\,. 
\end{array}\]
Again by the IFT, these equations will have a smooth solution branch in
the form $x=x(y)$ (that is, transverse to $S_0$ in $U$)
if and only if the Jacobian matrix of the map
\begin{equation}
S\to\br^d:x\mapsto(r_1,\ldots,r_{l-1},d_l,
            \ldots,d_{n-1},\del_ng_{n+1},\ldots,\del_ng_{d+1})
\end{equation}
(there are no $d_j$ terms if $n=l\,$) is nonsingular at $x_0\,$.  Now suppose
\begin{equation}  \label{e:rne0}
r_n:=r_n(x_0,0)\ne0
\end{equation}
which is satisfied automatically in the uniform case.  
Then solving $\til F_n(\eps,u)=0$
shows that $\del_n(x_0,0)\ne0$. Differentiating $\del_ng_i$ at
$(0,x_0,0)$ gives
\[
  \partial_x(\del_n g_i)(0,x_0,0) = \del_n(x_0,0) \partial_xg_i(0,x_0,0)
\]
for $i>n$ as $g_i(x_0,0)=0$.  Therefore, since $\del_n(x_0,0)\ne0$,
the Jacobian condition above can be re-stated as follows.
\begin{prop}
\label{p:sufficient}
Given that the necessary conditions (\ref{e:nconds})
are satisfied at $x_0\in S$ and~(\ref{e:rne0}) also holds, 
a sufficient condition for the existence of a solution
branch emanating from $(x_0,0)$ and transverse to $S_0$ in $K$ is that 
$x_0$ be a regular zero of the map
\[
g_{l,n}:=(r_1,\ldots,r_{l-1},d_l,\ldots,d_{n-1},g_{n+1},\ldots,g_{d+1}):
                  S\to\br^d,
\]
the $r$ and $d$ terms being absent when $n=1$.
\end{prop}
In the uniform case (where all $r_j\ne0$) the expression is again much simpler:
since $m_1=\cdots=m_n$ and there are no $r_j$ terms we have:
\begin{cor}  \label{c:suff}
In the uniform case, a sufficient
condition for $x_0\in S$ to be a branch point is that $x_0$ be a
regular zero of the map $g_{[n]}:S\to\br^d$.
\end{cor}
As before, note that if $m_1<m_2$ then $g_{[n]}$ 
becomes the map $\bar g:(g_2,\ldots,g_{d+1}):S\to\br^d$.
\medskip

\rem  
We can regard the normally {\em nondegenerate} case as an extension of
the above: here $m_1=1$ and all other $m_i>1$, and so regular
zeros of $\bar g$ are branch points.
\medskip

\rem
In general if a zero $x_0$ of $g_{l,n}$ is not a regular zero then
progress can still be made by standard bifurcation techniques.
For example, if the Jacobian matrix at $x_0$ has rank $d-1$ 
then a further Lyapunov-Schmidt reduction  
leads to a single \emph{bifurcation function} of one variable, 
whose zero set can be determined locally from the quadratic terms if
nondegenerate or from higher order terms by standard
singularity-theory methods otherwise.
\subsection{Examples}  \label{ss:examples}
We illustrate the preceding formal descriptions by some simple examples. 
\medskip

\begin{description}
\item[Example 1]\quad($\dim S=1$) \qquad $\til F_0(x,y)=(r_1(x,y)y^2,r_2(x,y)y^3)$.
\bigskip

The equations to solve for the zero set of the deformation $\til F_\eps(x,y)$ are
\begin{equation}  
\label{e:ex1} 
\begin{array}{ll}
  \eps g_1(\eps,x,y)+r_1(x,y)y^2 =0 \\[3mm]
  \eps g_2(\eps,x,y)+r_2(x,y)y^3 =0.
\end{array}
\end{equation}
From Proposition~\ref{p:necessary}
a necessary condition for $x_0\in S$ to be a branch point
is that $x_0$ be a zero of the function
\[
x\mapsto r_1(x,0)g_2(0,x,0)
\]
and a sufficient condition is that $x_0$ be a simple zero.  
\medskip

\item[Example 2]\quad($\dim S=1$) \qquad $\til F_0(x,y)=(r_1(x,y)y^2,r_2(x,y)y^2)$.
\bigskip

With $(g_1,g_2)$ as in (\ref{e:ex1}) a necessary condition now is that
$x_0$ be a zero of the function
\[
d_{12}:x\mapsto g_1(0,x,0)r_2(x,0)-g_2(0,x,0)r_1(x,0)
\]
and a sufficient condition is that $x_0$ be a simple zero.
\end{description}
\medskip

In the uniform case the functions $r_1,r_2$ are nonzero constants
and (as we have seen in general)
the necessary and sufficient conditions are conditions on $g_1,g_2$ alone.
Generically $g_2$ (Example 1) or $d_{12}$ (Example 2)
will have a discrete set of zeros on $S\,$, all of them simple.
\medskip

\begin{description}
\item[Example 3]\quad($\dim S=2$) \qquad  
           $\til  F_0(x,y)=(r_1(x,y)y^2,r_2(x,y)y^3,r_3(x,y)y^4)$.
\bigskip

Here with the usual notation we see that necessary conditions for 
$x_*\in S$ to be a branch point are that $r_1g_2=r_1g_3=r_2g_3=0$ at
$(x_*,0)\in S\times\br$, or in other words that $x_*$ is a zero either
of $(r_1g_2,g_3)$ or of $(r_1,r_2)$, these maps $S\to\br^2$ being evaluated with
$y=0$.  In the uniform case this reduces to $x_*$ being a zero
of $(g_2,g_3)$; then being a regular zero is a sufficient
condition for $x_*$ to be a branch point.
\end{description}

We leave to the reader the extension of this example to the cases when
some components of $\til F_0$ have leading terms with equal powers of $y\,$.
\subsection{The variational case}  \label{ss:variationalcase2}
We now consider how this local branching analysis applies in the
variational case.  For simplicity we take $q=k=1$.   
\medskip

Following the preceding sections
and dropping the tilde notation let
$f:\br\times U \to\br$ be a smooth function such
that $df_0(u)=0$ and $d^2f_0(u)=0$ for all $u=(x,0)\in S_0$,
where $f_\eps(u)=f(\eps,u)$.
Assuming the non-flatness condition (A5) we write $f$ in the form
\[
f(\eps,x,y)=\eps g(\eps,x,y) + y^m r(x,y)
\]
with $m\ge3$ and $r(\,\cdot\,,0)$ not identically zero on $S\,$;
in the uniform case $r({\cdot},0)\ne0$ on $S$.  
Critical points of $f_\eps:U\to\br$ occur where
\begin{equation}
\label{e:varEQ}
\left\{\begin{array}{ll}
\eps\,\partial_xg(\eps,x,y)+y^m\,\partial_xr(x,y)=0\in \br^d\\[3mm]
\eps\,\partial_yg(\eps,x,y)+m\,y^{m-1}r(x,y)+y^m\,\partial_yr(x,y)=0\in\br.
\end{array}\right.
\end{equation}
\begin{prop} \label{p:necvar}
A necessary condition for $x_0\in S$ to be a branch point is
\begin{equation}
r(x_0,0)\,\partial_xg(0,x_0,0)=0\in\br^d.
\label{e:varNC}
\end{equation}
\end{prop}
\proof
The argument of Proposition~\ref{p:necessary} leads (after
cancelling $y^{m-1}$) to the equation
\[
m\,r(x,y)\partial_xg(\eps,x,y)+y\left[\partial_xg(\eps,x,y)\,
\partial_yr(x,y)-\partial_yg(\eps,x,y)\,\partial_xr(x,y)\right]=0
\]
and as $(\eps,x,y)\to(0,x_0,0)$ we obtain~(\ref{e:varNC}).
\endproof
\begin{cor}
In the uniform case all branch points are critical points of
$g(0,{\cdot},0): S_0\to\br$. \endproof
\end{cor}
\begin{prop}  \label{p:varbranch}
If (\ref{e:varNC}) holds for $x_0\in S$, then 
a sufficient condition for $x_0$ to be a branch point is provided 
by either of the following conditions at $x_0\,$:
\begin{enumerate}
\item
$\partial r/\partial x_{n-1}\ne0$ and $x_0$ is a regular zero of the map
$$x \mapsto (mr,e_2,\ldots,e_{n-1},h_n,\ldots,h_d):S\to\br^d$$
where $e_j:=-\dbd gy\dbd r{x_{j-1}}$ and $h_i=\dbd g{x_i}$;
\item
$r\ne0$ and $x_0$ is a nondegenerate critical point of $g(0,{\cdot},0)\,$.
\end{enumerate}
\end{prop}
\proof
If $\partial_xg(0,x_0,0)\ne0$ then $r(x_0,0)=0$ by~(\ref{e:varNC}).
Condition~(\ref{e:rne0}) here is $\partial r/\partial x_{n-1}\ne0$,
where (with the notation of Section~\ref{ss:nec}) we have $J_n=[2,d+1]$
and (i) follows as in Proposition~\ref{p:sufficient}: note that $d_j$
reduces to $e_j$ since $r=0$.

If $\partial_xg(0,x_0,0)=0$ then again from
Proposition~\ref{p:sufficient} with~(\ref{e:rne0}) we obtain (ii).
\endproof
\begin{cor}
In the uniform case $r\ne0$ and it is only the condition
(ii) that applies.  \endproof
\end{cor}
Here, just as in the normally nondegenerate case,
we can use topological tools such as 
Lyusternik-Schnirelmann category or Morse Theory
in order to give a lower bounds for the
number of branch points: see~\cite{ACE,Chang,RE2}.
Observe that this global result characterising the branch
points is the same as in the normally nondegenerate case.
However, the branching behaviour itself is different: 
recall the Remark following Theorem~\ref{t:bif2}.
\section{Applications}  \label{s:apps}
Here we outline some examples of problems where the techniques presented 
in this paper can be applied. 
%
\subsection{Periodic orbits}
A standard technique for finding periodic orbits (of period $T$) for a
vector field $X_0$ is to consider them as zeros of the nonlinear
operator 
\[
F_0:=\frac{d}{dt}-X_0
\] 
defined on a suitable Banach space $\mathcal H$ of $T$-periodic functions.
There is a natural action of the circle group $\bs^1$ on $\mathcal H$
given by time-translation modulo $T$, and if $F_0(x)=0$
then also $F_0(x_s)=0$ for all $s\in \bs^1=\br\bmod T$, where 
$x_s(t)$ denotes $x(s+t)$.  Thus each zero of $F_0$ is
automatically a member of an $\bs^1$-orbit $S$ of zeros.  
When this orbit is normally nondegenerate the problem essentially
reduces to a problem on $S$ and techniques that are now standard
(\cite{ACE,C4,CC3,D2,HA2,V1}) can be applied to study
bifurcations of zeros of $F_0$ under perturbations of $X_0\,$.   
The first authors apparently to attempt a systematic
approach to normal degeneracy in this context were Hale and
Taboas \cite{HT2} who used analytic methods to study a differential
equation of the form
\[
\ddot x + g(x) = \eps_1h(t)\dot x+\eps_2f(t) \in \br
\]
close to a degenerate $T$-periodic solution of $\ddot x+g(x)=0$, where
$h$ and $f$ are $T$-periodic functions.
As indicated in Section~\ref{sss:bifgeom} our
Theorem~\ref{t:bif1} recovers their results and sets them in a
geometric context. 
\subsection{Hamiltonian systems}
In the case of
Hamiltonian systems, where periodic orbits may typically lie in
flow-invariant tori, each zero of $F_0$ will lie in a
larger-dimensional manifold (torus) of zeros.  We expect 
normal non-degeneracy of this manifold within each energy level, 
but for discrete values of the Hamiltonian function (energy) we may
find normal degeneracy. This appears both in the study of the 
Poincar\'e map for the periodic orbits and equivalently
in the study of manifolds of critical points of the Jacobi metric.
\msk

To illustrate this in a relatively simple situation where the periodic
orbits are isolated in each energy level, consider a Hamiltonian
system in $\br^2$ with a smooth Hamiltonian of the form 
\[
H(q,p)=K(p)+V(q)
\]
where $(q,p)\in T^*\br^2\cong\br^2\times\br^2$, and where
$K(p)=|p|^2/2$.  Up to parametrisation, a  periodic orbit $q({\cdot})$ 
for the system is a critical point of the Jacobi metric 
(see~\cite{AM} for example) 
\[
\J[q]=\int_0^1\,(E-V(q(s)))\,K(\dot q(s))ds
\]
which has natural time-translation $\bs^1$-invariance as a smooth
function on ${\cal H}=H^1([0,1],\br^2)$;
here $E$ is a constant greater than $V(q)$ throughout the region of
$\br^2$ under consideration.
\medskip

Assume that the system is rotationally symmetric about the
origin, so that $V(q)=V(r)$ with $r=|q|$.  The Euler-Lagrange equation
derived from $\J$ is 
\begin{equation} \label{e:lag}
 \frac{d}{dt}[\dot q\,(E-V(r)] 
    + \frac{|\dot q|^2}{2}\, \frac{1}{r}V'(r)\,q(t)=0
\end{equation}
where ${'}$ denotes $\dbd{}r$.
Now let $q_0(t)$ be a uniform circular orbit written as $q_0(t)=r_0\,e(t)$
where $e(t)\in\br^2$ is a unit vector.
Substituting $q_0(t)$ into~(\ref{e:lag}) and using
the identity $\ddot e(t)=-|\dot e(t)|^2\,e(t)$ gives
\begin{equation}  \label{e:circular}
(E-V(r_0))=\frac{r_0}{2}V'(r_0)
\end{equation}
as a necessary and sufficient condition for $q_0(t)$ to be a periodic
orbit for the system.
Note that circular periodic orbits of fixed radius form a manifold $S_E$
in $\mathcal H$ that is a copy of $\bs^1$:
\[
S_E=\{q_s\in{\mathcal H}: q_s(t)=R_s\,q_0(t)\}
\] 
where $R_s\in SO(2)$ denotes rotation through angle $s\in [0,2\pi)$.

To ascertain normal (non)degeneracy we next study the second derivative
$D^2\J[q_0]$ which is defined on on the space $\mathcal H$ and
corresponds to a densely defined self-adjoint operator on $L^2([0,1],\br^2)$.
We show in the Appendix (Section~\ref{ss:ham}) that given~(\ref{e:circular}) 
the kernel of $D^2\J[q_0]$ has dimension $2$ when
\begin{equation}    \label{e:condition}
r_0\,V''(r_0) = (2n^2-9)V'(r_0)
\end{equation}
for some $n\in \bz$, and has dimension $1$ (corresponding to $T_{q_0}S$) 
otherwise. 
Note that condition (\ref{e:condition}) can be read in terms of
discrete values $\{E_n\}$ of the energy $E\,$:
\begin{equation}
E_n=V(r_0)+\frac{r_0^2}{2(2n^2-9)}\,V''(r_0).
\label{e:energy}
\end{equation}
For energy values $E_n$ the corresponding $S_{E_n}$ is 
a normally degenerate manifold of points representing periodic orbits.
\msk

A potential $V(q)$ for which all this can be explicitly verified is
the so-called {\em Mexican hat}:
\[V(q)=-\frac{\lambda^2}{2}|q|^2+\frac{1}{4}|q|^4.\]
The above methods allow us to detect the presence of periodic orbits
that persist after perturbing $V$
in the annular \emph{Hill's region} given by 
$\{q\in\br^2:V(q)\le E\}$ where $-\frac14\lambda^4<E<0$. 
Setting $E=E_n$  and fixing $u\in K_{q_0}$ where $K_{q_0}$ is a complement to
$T_{q_0}S_{E_n}$ in $\ker D^2\J[q_0]$ it is not difficult to check that the Taylor
expansion about $q_0$ of the Jacobi metric after Lyapunov-Schmidt
reduction reads 
\begin{equation}    \label{eq:reduced-jacobi}
\J(y):=\J[q_0+y\,u]=\J[q_0]+\J_0 y^4 + O(y^5)
\end{equation}
where $\J_0=\frac1{4!}D^4\J[q_0](u,u,u,u)$.  This places us in the context of
Section~\ref{ss:varcase} with $d=k=1$ and $m=4$.

An example of bifurcation can be seen by considering a potential
\[V_\eps(q)=V(q)+\eps_1\phi_1(q)+\eps_2\phi_2(q)\]
where $\phi_i(q)=\phi_i(q_1,q_2)$ is a smooth function, $i=1,2$.
Up to degree 4 in $y$ the reduced Jacobi metric is now
\begin{equation}  \label{eq:reduced-jacobi2}
\J_\eps[q_0+y\,u]=\J[q_0]+\J_0y^4+\eps_1\Phi_1[q_0,y]+\eps_2\Phi_2[q_0,y]
\end{equation}
where 
\[
\Phi_i[q_0,y]=-\int_0^1 dt\|\dot{q}_0(t)\|^2(\phi_i(q_0(t)+y u(t))
\]
for $i=1,2$.  Now let $q_0(t)=(r_0\cos(\omega\,t+x),r_0\sin(\omega\,t+x))$ 
be an element of $S_{E_n}$ parametrised by $x\in S^1$.
The expression (~\ref{eq:reduced-jacobi2}) reduces to 
\[
\J(\eps,x,y):=\J[q_0]+\J_0\,y^4+\eps_1g_1(x,y)+\eps_2g_2(x,y)
\]
up to degree $4$ in $y$, where
\[g_i(x,y)=\Phi_i[q_0,y]\]
for $i=1,2$.  Following the Example given in Section~\ref{ss:varcase} but
now with $m=4$ we see that the bifurcation geometry is
determined by the geometry of the curves $B_0'$ and $B_1$ in 
$\bs^1\times S$ given by $b_0'(s,x)=0$ and $b_1(s,x)=0$ respectively,
where $b_0,b_1$ are derived from the perturbation $g$ as in
Section~\ref{ss:m>2} with $m=4$.
\msk

\rem The case of a single parameter $\eps$ (so $q=1$) can be
interpreted in this context as a fixed choice of $s\in\bs^1$.
Here the function $g(0,{\cdot},0):S_0\to\br$ of
Section~\ref{ss:variationalcase2} corresponds to the function
$b_0(s,{\cdot})$ of Sections~\ref{ss:varcase} and~\ref{ss:m>2}.
The necessary condition for $x_0$ to be a branch point as given by
Proposition~\ref{p:necvar} corresponds to the statement that
the curve $B_0'$ intersects the circle $\{s\}\times S$ at
$(s,x_0)$, and the sufficient condition given by
Proposition~\ref{p:varbranch}(ii) corresponds to the statement that
it does so transversely.  Since every smooth function on a circle has
at least two critical points there will be at least two such
intersections (Lyusternik-Schnirelmann category), and in any case an
even number if they are transverse (Morse theory).  
\subsection{Steady states in chemical reaction networks} \label{ss:chem}
In kinetic models for chemical reaction networks 
the time evolution of a vector $x(t)=(x_1(t),...,x_n(t))$ of 
concentrations of $n$ chemical species involved in $r$ reactions is given by
\begin{equation}
\frac{d x}{d t}=F(x,k):=B\,\nu(x,k)
\label{eq:reaction-eqs}
\end{equation}
where $B$ is an $n\times r$ matrix (\emph{stochiometric coefficients}) 
and 
\[
\nu:\br^n\times\br^c\to\br^r
\]
is a smooth map, often in fact a polynomial map.
A review of this type of system can be found in \cite{Domijan}. 
The components of $k\in\br^c$ are the \emph{reaction constants}.  For
given $k$ the stationary states of (\ref{eq:reaction-eqs}) are thus
given by the set
\[S_k:=\{x\in\br^n:\nu_k(x)\in\ker(B)\}\]
where $\nu_k(x)=\nu(x,k)$.
\msk

In applications, one is interested in studying steady states of
perturbed systems of the form
\begin{equation}  \label{e:perchem}
F_\eps(x):=B\,\nu(x,k)+\eps\phi(\eps,x)
\end{equation}
for $\eps\in\br^q$ where $\phi(\eps,x)$ is smooth (not necessarily
polynomial) function.  It is therefore appropriate to apply our
analysis to systems of this type. 
\msk

Let $R\subset\br^n$ denote the range of $B$ and let $\pi:\br^n\to R$
be a projection onto $R$.  Write $\til B=\pi\circ B:\br^r\to R$.
If $0\in R$ is a regular value of $\til B\circ\nu_k$ (that is, 
the map $\nu_k$ is transverse to $\ker B$)
then $S_k$ is a smooth submanifold of $\br^n$
with codimension equal to the dimension of $R$ (the rank of $B$);
moreover, $\dim S_k = \dim\ker DF(x)$ for $x\in S_k$ so the manifold
$S_k$ is normally nondegenerate.
\msk

If $0$ is not a regular value of $\til B\circ\nu_k$ it means that there is
at least one point $x\in S_k$ at which $\nu_k$ fails to be transverse
to $\ker B$.  Situations can arise where this occurs simultaneously at
all points of $S_k$ and we have uniform normal degeneracy.
We present now a simple example of this class of system.
\msk

Consider a set of reactions of the form
 \[
 X_1+X_2\mathop{\rightarrow}\limits^{1} X_3\,,
  ~~X_3 \mathop{\rightarrow}\limits^{1} X_1+X_2\,,
  ~~\emptyset\mathop{\longrightarrow}\limits^{v([X_1])} X_2
\]
 where the first two reaction rates are fixed to $1$ and the third is
 a function $v$ of the concentration $[X_1]$ of $X_1$. 
 The third reaction is interpreted 
 as an external input whose rate is a function of $[X_1]$. 
 Now let $[X_i]=x_i$  denote the $i$-th concentration, 
 so  $x=(x_1,x_2,x_3)\in\br_+^3$ (all
 components non-negative) and the dynamic equations read 
 \begin{equation}
 \begin{array}{lll}
 \dot{x}_1(t)=-x_1(t)\,x_2(t)+x_3(t)\\[3mm]
 \dot{x}_2(t)=-x_1(t)\,x_2(t)+x_3(t) +v(x_1(t))\\[3mm]
 \dot{x}_3(t)=x_1(t)\,x_2(t)-x_3(t)
 \end{array}
 \label{e:reactions}
 \end{equation}
where here we take $k$ as fixed and omit it from the notation.  
The map $\nu:\br^3\to\br^3$ is given by
 \[\nu(x)=(x_1x_2\,,\,v(x_1)\,,x_3)^T,\]
and the matrix $B$ is 
 \[B=\left(\begin{array}{ccc}
 -1 & 0 & 1\\
 -1 & 1 & 1\\
 1 & 0 & -1
 \end{array}\right)\]
with $1$-dimensional kernel $\,\ker B=\spn\{(1,0,1)^T\}\,$ and range
$\,R=\spn\{(0,1,0)^T,(1,1,-1)^T\}$.
\msk

 Assume that the rate $v({\cdot})$ is a smooth function with 
 a unique zero $x_1^*>0$. The stationary states are therefore
\[
S_+=\{(x_1,x_2,x_3): x_1x_2=\lambda,~~v(x_1)=0,
                  ~~x_3=\lambda\}
\]
for $\lambda>0$, which is the part of the $1$-dimensional affine space
 \[S=\{(x_1,x_2,x_3): x^*_1x_2-x_3=0,~~x_1=x^*_1\}\]
where all components are positive.  The tangent space of $S$ at $x\in
S$ is independent of $x$, namely
\[
T_xS=\spn\{(0,1,x_1^*)^T\}.   
\]
We have
\[
D\nu(x)=\left(\begin{array}{ccc}
 x_2 & x_1 & 0\\
 v'(x_1) & 0 & 0\\
 0 & 0 & 1
 \end{array}\right)
\]
 and
\[
BD\nu(x)=\left(\begin{array}{ccc}
 -x_2 & -x_1& 1\\
 v'(x_1)-x_2 & -x_1 & 1\\
 x_2 & x_1 & -1
 \end{array}\right).
\]
 Hence we see
\begin{prop}  \label{p:reg}
The map $\til B\nu:\br^3\to R$ has $0$ as a regular value if and only
if $v'(x_1^*)\ne0$.
\endproof
\end{prop}
Now suppose that $v'(x_1^*)=0$.  This does not alter $S$ but does
destroy the transversality of $\nu$ to $\ker B$.  
If $x\in S_+\subset S$ then
$T_xS$ is annihilated by $BD\nu(x)$ but not by $D\nu(x)$,
 so that $S_+$ is normally degenerate with constant corank $1$.
We choose a suitable complement $K_x$ to $T_xS$ in $\ker BD\nu(x)$.
In the Appendix Section \ref{sect:degeneracy} it is shown that the 
reduced equation $\til F_0(u)=0$  (see Section~\ref{sect:LS})
has uniform quadratic degeneracy provided $v''(x_1^*)\ne0$.  For studying
bifurcations we are therefore in the context of Section~\ref{ss:m=2}
with $d=k=1$ and $m=2$.
\section{Appendix: calculations for examples in Section~\ref{s:apps}.} 
\subsection{Hamiltonian periodic orbits}  \label{ss:ham}
 In this section we derive the condition (\ref{e:condition}).
\msk

Let $q(t)=r\,e(t)\in\br^2$ be a circular path with constant speed 
$|\dot e(t)|=\omega$.  Let $u(t),v(t)\in{\mathcal H}$.
Direct calculation shows that the Hessian of $\J$ at $q$ is

\[\begin{array}{lll}
\displaystyle D^2\J[q](u,v)=\int_0^1dt\left[\langle \dot u(t),\dot
  v(t)\rangle\,(E-V(r))\right]+\\[5mm]
\displaystyle-\,\omega\int_0^1dt\left[\langle
  v(t),e(t)\rangle\langle\dot u(t),e^\perp(t)\rangle
  +\langle u(t),e(t)\rangle\langle\dot
  v(t),e^\perp(t)\rangle\right]rV'(r)+\\[5mm]
\displaystyle
  -\frac{\omega^2}2\int_0^1dt\left[r^2\left(V''(r)+
\frac1r V'(r)\right)\langle e(t),u(t)\rangle\,\langle e(t),v(t)\rangle 
+rV'(r)\langle u(t),v(t)\rangle\right].
\end{array}\]

\noindent To compute $D^2\J$ at $q=q_0$ we set $r=r_0$ and use 
condition (\ref{e:circular}).  We find
\[\begin{array}{lll}
\displaystyle D^2\J[q_0](u,v)=\int_0^1dt\left[\langle \dot
  u(t),\dot v(t)\rangle\, (E-V(r_0))\right]+\\[5mm]
\displaystyle -2\omega\int_0^1dt\left[\bigl(\langle v(t),e(t)\rangle
\langle\dot u(t),e^\perp(t)\rangle + \langle u(t),e(t)\rangle
\langle\dot v(t),e^\perp(t)\rangle\bigr)\bigl(E-V(r_0)\bigr)\right]+\\[5mm]
\displaystyle -\omega^2\int_0^1dt\left[\frac12\left(r_0^2V''(r_0)
  +2\,(E-V(r_0))\right)\langle e(t),u(t)\rangle\,
   \langle e(t),v(t)\rangle+\bigl(E-V(r_0)\bigr)\,\langle v(t),u(t)\rangle\right].
\end{array}\]
Since $E-V(r_0)>0$ the condition $D^2\J[q_0](u,v)=0$ is equivalent to
\begin{equation}
\begin{array}{lll}
\fl\displaystyle \int_0^1dt\left[\langle \dot v(t),\dot u(t)\rangle
-2\omega\bigl(\langle v(t),e(t)\rangle\langle\dot
u(t),e^\perp(t)\rangle+\langle u(t),e(t)\rangle
\langle\dot v(t),e^\perp(t)\rangle\bigr)+\right.\\[5mm]
\qquad\displaystyle \left.\qquad -(\omega^2/2)\,\Omega\,\langle e(t),u(t)\rangle\,
\langle e(t),v(t)\rangle 
-\omega^2\,\langle v(t),u(t)\rangle\right]=0
\end{array}\label{e:D2}
\end{equation}
where
\begin{equation}   \label{e:defthisK}
\Omega:= \left(\frac{r_0^2V''(r_0)}{2\,(E-V(r_0))}\,+1\right).
\end{equation}

Integration by parts and a standard argument from calculus of variations
shows that $u(t)$ lies in the kernel of $D^2\J[q_0]$ precisely when
\begin{equation}
\begin{array}{lll}
\displaystyle \ddot u(t)
+2\omega\,e(t)\,\langle\dot u(t),e^\perp(t)\rangle
     -2\omega\,\frac{d}{dt}(\langle u(t),e(t)\rangle)\,e^\perp(t)+\\[5mm]
\displaystyle\qquad +2\,\omega^2\,\langle u(t),e(t)\rangle\,e(t)
     +\frac{\omega^2}{2}\,\Omega\,\langle e(t),u(t)\rangle\,e(t) 
+\omega^2\,u(t)=0.
\end{array}\label{e:D2-kernel-eq}
\end{equation}

Decomposing $u(t)$ as
\[u(t)=a(t)\,e(t)+b(t)\,e^\perp(t)\]
we find the equation (\ref{e:D2-kernel-eq}) splits into
\begin{equation}
\left\{\begin{array}{ll}
\ddot a(t)+ (4+\Omega/2)\,\omega^2\,a(t)=0,\\
\ddot b(t)=0.
\end{array}\right.
\label{e:final}
\end{equation}
Equations (\ref{e:final}) imply that
\begin{equation}
\dim \ker D^2\J[q_0]=\left\{\begin{array}{ll}
1\mbox{ if $(4+\Omega/2)\neq n^2$ \quad with \quad $n\in\bz$},\\
2\mbox{ if $(4+\Omega/2)= n^2$ \quad for some \quad $n\in\bz$}.\\
\end{array}\right.
\label{e:kernel}
\end{equation}
Finally, the expression~(\ref{e:defthisK}) for $\Omega$ and 
condition (\ref{e:circular}) show that the condition 
$(4+\Omega/2)= n^2$ is equivalent to~(\ref{e:condition}).
\endproof
\subsection{The reduced equation for the chemical reaction network}    
\label{sect:degeneracy}
Here we construct the reduced equation
for~(\ref{e:reactions}) perturbed as in~(\ref{e:perchem}). 
From~(\ref{e:reactions}) we write~(\ref{e:perchem}) as
\begin{equation}
 F_\eps(x):= B\nu(x) +\phi(\eps,x) = f(x)\,w + v(x_1)\,e+\eps\phi(\eps,x)
\end{equation}
where $w=(1,1,-1)^T\,$, $\,e=(0,1,0)^T$ and $f(x)=-x_1x_2+x_3\in\br$.
For $x\in S_+$ the range $R_x$ of $BD\nu(x)$ is $\spn\{w,e\}$ except if $v'(x_1^*)=0$
in which case $R_x=\spn\{w\}$. For $z\in\br_+^3$ we have 
\begin{equation}  \label{e:xyz}
\begin{array}{ll}
\langle w,F_\eps(z)\rangle &=f(z)\|w\|^2 
   + v(z_1)\langle w,e\rangle + \langle w,\phi(\eps,z)\rangle \\[3mm]
       &= 3f(z) + v(z_1) + \langle w,\phi(\eps,z)\rangle.
\end{array}
\end{equation}
For $x\in S_+$ and $v'(x_1^*)=0$ we have
$\ker BD\nu(x)=u^\perp$ where $u=(x_2,x_1,-1)^T$ and so it is
natural for the Lyapunov-Schmidt reduction to take $L_x=\spn\{u\}$.
Accordingly for $z$ in a \nhd of $S_+$ in $\br_+^3$ we write   
$$z=x+yn+\lambda u \in K_x\oplus L_x$$ 
where $x\in S_+$ and
$$K_x=(T_xS\oplus \spn\{u\}\bigr)^\perp=\spn \{n\}$$
where $n=(1+{x_1^*}^2,-sx_1^*,s)^T.$
By the IFT we can solve
\begin{equation}  \label{e:solvelambda}
  \langle w,F_\eps(x+yn+\lambda u)\rangle=0
\end{equation}
in the form $\lambda=\lambda^*(x,y,\eps)$ 
where $\lambda^*$ is a smooth function defined on a \nhd of 
$S_0\times 0$ in $\br^3\times\br^q$ and with $\dbd{\lambda^*}y(x,0,0)=0$ for
all $x\in S$. It then remains to solve the reduced equation pair
\begin{equation} \label{e:solvez}
\langle s,F_\eps(x+yn+\lambda^*(x,y,\eps)u)\rangle = 0,\qquad
\langle\/ p,F_\eps(x+yn+\lambda^*(x,y,\eps)u)\rangle = 0
\end{equation}
where $\spn\{s\}=T_xS$ and $\spn\{p\}=P_x$ chosen so that 
$P_x$ does not lie in $R_x\oplus T_xS$ (for example, $p=(0,1,0)^T$).  Now 
\[
\begin{array}{ll} \label{e:sF1}
\langle s,F_\eps(z)\rangle &= f(z) \langle s,w\rangle
             + v(z_1)\langle s,e\rangle +\eps\langle
             s,\phi(\eps,z)\rangle \\[3mm]
    &=a(s)v(z_1) + \eps \langle b(s),\phi(\eps,z)\rangle
\end{array}
\]
where 
\begin{equation}   \label{e:sF2}
    a(s)=\langle s,e\rangle -\frac13\langle s,w\rangle \qquad
\mbox{and} \qquad    b(s)=s-\frac13\langle s,w\rangle w
\end{equation}
from~(\ref{e:xyz}), and similarly for $\langle\/ p,F_\eps(z)\rangle$.  
With $z=x+yn+\lambda^*u$ we can write
\[
\begin{array}{ll}
v(z_1)=v(x_1+yn_1+\lambda^*u_1)&= v(x_1^*+y(1+{x_1^*}^2)+\lambda^*x_2)\\
          &= \frac12v''(x_1^*)(y(1+{x_1^*}^2)+\lambda^*u_1)^2 + O(3)
\end{array}
\]
where $O(3)$ terms are of third or higher order in $(y_1,\eps)$ and
smoothly parametrised by $x\in S_+$.  Since $\lambda^*$ has no linear
term in $y$ the smooth map $y_1\mapsto y(1+{x_1^*}^2)+\lambda^* u_1$ is a
diffeomorphism germ at the origin in $K_x=\br$ smoothly parametrised by
$x\in S_+$ and sufficiently small $\eps\in\br^q$.  After we make this
change of coordinate, the assumption $v''(x_1^*)\ne0$ allows us to apply
and a further local diffeomorphism in $K_x$ to absorb the $O(3)$ terms. We may
therefore suppose without loss of generality that the equations~(\ref{e:solvez})
are
\[
\begin{array}{ll}
  a(s)cy_1^2 + \eps \langle b(s),\phi(\eps,x,y)\rangle \\
  a(p)cy_1^2 + \eps \langle b(p),\phi(\eps,x,y)\rangle 
\end{array}
\]
where $c=\frac12 v''(x_1^*)\ne0$, and so we have a quadratic
system as claimed in Section~\ref{ss:chem}.  More specifically, the
problem is of the form given in Example 2 in Section~\ref{ss:examples}
and hence the branch points $x\in S_+$ are given by the
zeros (provided they are simple) of the function
\[
  g(x):= a(s)\langle b(p),\phi(0,x,0)\rangle
       - b(s)\langle a(p),\phi(0,x,0)\rangle.
\]
\bsk

{\em Acknowledgements.}
The authors are grateful for partial support from the London
Mathematical Society through a Scheme 4 grant, ref. 4433.  We also
indebted to David Kirby for helpful conversations about the algebra in
Section~\ref{ss:m>2}.


\end{document}